# ON THE OVERLAP IN THE MULTIPLE SPHERICAL SK MODELS[1]

By Dmitry Panchenko and Michel Talagrand

*Massachusetts Institute of Technology, Universite Paris 6 and Ohio State University*

In order to study certain questions concerning the distribution of the overlap in Sherrington–Kirkpatrick type models, such as the chaos and ultrametricity problems, it seems natural to study the free energy of multiple systems with constrained overlaps. One can write analogues of Guerra's replica symmetry breaking bound for such systems but it is not at all obvious how to choose informative functional order parameters in these bounds. We were able to make some progress for spherical pure $p$-spin SK models where many computations can be made explicitly. For pure 2-spin model we prove ultrametricity and chaos in an external field. For the pure $p$-spin model for even $p > 4$ without an external field we describe two possible values of the overlap of two systems at different temperatures. We also prove a somewhat unexpected result which shows that in the 2-spin model the support of the joint overlap distribution is not always witnessed at the level of the free energy and, for example, ultrametricity holds only in a weak sense.

**1. Introduction and main results.** Let us consider a Gaussian–Hamiltonian (process) $H_N(\boldsymbol{\sigma})$ indexed by $\boldsymbol{\sigma} \in \mathbb{R}^N$ with covariance that satisfies

$$\left| \frac{1}{N} \mathbb{E} H_N(\boldsymbol{\sigma}^1) H_N(\boldsymbol{\sigma}^2) - \xi(R_{1,2}) \right| \leq c_N, \tag{1.1}$$

where $c_N \to 0$, $R_{1,2} = N^{-1} \sum_{i \leq N} \sigma_i^1 \sigma_i^2$ is the overlap of configurations $\boldsymbol{\sigma}^1, \boldsymbol{\sigma}^2$ and $\xi$ is a smooth enough convex even function with $\xi(0) = 0$. We define $\theta(x) = x\xi'(x) - \xi(x)$. Even though we will state some basic results for a general function $\xi$, our main results will deal with the pure $p$-spin SK Hamiltonians that correspond to $\xi(q) = q^p/p$ for even $p \geq 2$. For example, one can

Received April 2006; revised November 2006.

[1]Supported in part by NSF grants.

consider

(1.2) $$H_N(\boldsymbol{\sigma}) = \frac{1}{\sqrt{p}N^{(p-1)/2}} \sum_{1 \leq i_1,\ldots,i_p \leq N} g_{i_1,\ldots,i_p}\sigma_{i_1}\cdots\sigma_{i_p},$$

where $(g_{i_1,\ldots,i_p})$ are i.i.d. Gaussian random variables. For $p=2$ this is a classical SK Hamiltonian [5]. The factor $p^{-1/2}$ is not important and is chosen so that $\xi'(q) = q^{p-1}$. In this paper we will consider the spherical model when the spin configurations $\boldsymbol{\sigma}$ belong to a sphere $S_N$ of radius $\sqrt{N}$ and the a priori distribution of $\boldsymbol{\sigma}$ is the uniform measure $\lambda_N$ on $S_N$. Given an inverse temperature $\beta > 0$ and an external field $h \in \mathbb{R}$, the "free energy" is defined by

$$F_N(\beta,h) = \frac{1}{N}\mathbb{E}\log\int_{S_N} \exp\left(\beta H_N(\boldsymbol{\sigma}) + h\sum_{i\leq N}\sigma_i\right)d\lambda_N(\boldsymbol{\sigma}).$$

Its limit $\lim_{N\to\infty} F_N(\beta,h) = \mathcal{P}(\beta,h)$ was computed in [1] and the computation was made rigorous in [8]. The limit $\mathcal{P}(\beta,h)$ is the analogue of the Parisi formula in the Ising SK model [6, 7]. For the spherical model it is somewhat simplified by the fact that certain computations become more explicit. Namely, if, given $k \geq 1$ and given two sequences $\mathbf{m} = (m_l)_{0\leq l\leq k}$ and $\mathbf{q} = (q_l)_{0\leq l\leq k+1}$ such that

$$0 = m_0 \leq m_1 \leq \cdots \leq m_k = 1,$$
$$q = q_0 \leq q_1 \leq \cdots \leq q_{k+1} = 1,$$

and given a parameter $b > 1$, we define for $l \leq k$

$$d_l = \sum_{l\leq p\leq k} m_p(\xi'(q_{p+1}) - \xi'(q_p)) \quad \text{and} \quad D_l = b - d_l,$$

then

(1.3) $$\mathcal{P}(\beta,h) = \inf_{b,k,\mathbf{m},\mathbf{q}} \frac{1}{2}\left(b - 1 - \log b + \frac{1}{D_1}(h^2 + \xi'(q_1))\right.$$
$$\left. + \sum_{1\leq l\leq k} \frac{1}{m_l}\log\frac{D_{l+1}}{D_l} - \sum_{1\leq l\leq k} m_l(\theta(q_{l+1}) - \theta(q_l))\right).$$

$\mathcal{P}(\beta,h)$ is also given by the Crisanti–Sommers representation (see Section 4 in [8]) as follows. If $\delta_l = \sum_{l\leq p\leq k} m_p(q_{p+1} - q_p)$ then

(1.4) $$\mathcal{P}(\beta,h) = \inf_{k,\mathbf{m},\mathbf{q}} \frac{1}{2}\left(h^2\delta_1 + \frac{1}{\delta_1}q_1 + \sum_{1\leq l\leq k-1} \frac{1}{m_l}\log\frac{\delta_l}{\delta_{l+1}} + \log\delta_k\right.$$
$$\left. + \sum_{1\leq l\leq k} m_l(\xi(q_{l+1}) - \xi(q_l))\right).$$

When $h = 0$ we will write $\mathcal{P}(\beta) := \mathcal{P}(\beta, 0)$. The goal of this paper is to prove and analyze some bounds on the free energy of multiple copies of the system, possibly at different temperatures, coupled by constraining their overlap. Let $Q$ be a $n \times n$ symmetric nonnegative definite matrix with elements $q_{j,j'} \in [-1, 1]$ and $q_{j,j} = 1$. Given $\varepsilon > 0$ consider a set

$$(1.5) \quad Q_\varepsilon = \{(\boldsymbol{\sigma}^1, \ldots, \boldsymbol{\sigma}^n) \in S_N^n : R_{j,j'} \in [q_{j,j'} - \varepsilon, q_{j,j'} + \varepsilon] \text{ for } j, j' \leq n\}$$

and given $\beta_1, \ldots, \beta_n > 0$ and $h_1, \ldots, h_n \in \mathbb{R}$ we define a "free energy" of the $n$-configuration system constrained to the set $Q_\varepsilon$ by

$$(1.6) \quad F_N(Q_\varepsilon) = \frac{1}{N} \mathbb{E} \log \int_{Q_\varepsilon} \exp\left( \sum_{j \leq n} \beta_j H_N(\boldsymbol{\sigma}^j) + \sum_{j \leq n} h_j \sum_{i \leq N} \sigma_i^j \right) d\lambda_N^n.$$

Obviously, $F_N(Q_\varepsilon) \leq \sum_{j \leq n} F_N(\beta_j, h_j)$ and as a result a trivial bound would be

$$(1.7) \quad \limsup_{N \to \infty} F_N(Q_\varepsilon) \leq \sum_{j \leq n} \mathcal{P}(\beta_j, h_j).$$

We would like to construct some nontrivial bounds on $F_N(Q_\varepsilon)$ that would yield some information on the support of the distribution of the overlaps $(R_{j,j'})$ under the product Gibbs measure by showing, for example, that for some constraints $Q$

$$\limsup_{\varepsilon \to 0} \limsup_{N \to \infty} F_N(Q_\varepsilon) < \sum_{j \leq n} \mathcal{P}(\beta_j, h_j),$$

which by concentration of measure would imply that with high probability the overlaps cannot be in configuration $Q$ for the product Gibbs measure. At this moment, the only approach we could conceive for proving such bounds is based on an analogue of Guerra's interpolation [3] that was used in [7] for two coupled systems at equal temperatures and external fields, and explained in more generality in [9]. It was also explained in [9] that it seems to be not obvious at all how to choose parameters in these bounds that would at least witness the obvious inequality (1.7). In this paper we will describe several situations when we were able to find such parameters. The methods of the proofs are at least as interesting as the results they imply since they shed some light on the difficulties of finding informative parameters in the bounds and give some hope that, in principle, these bounds might be "correct" and it could be only a (very difficult) technical problem to find suitable parameters in more general situations.

*The pure* 2-*spin SK model.* The first case we will consider is the pure 2-spin SK model in (1.2) with $p=2$ without external field, that is, $h=0$. What makes this case particularly simple is that due to the proof of Proposition 2.2 in [8] the infimum in (1.4) [and (1.3)] is achieved on the replica-symmetric choice of parameters, that is, for $k=1$, so that (1.4) becomes

$$\mathcal{P}(\beta) = \inf_{q\in[0,1]} \frac{1}{2}\left(\frac{q}{1-q} + \log(1-q) + \beta^2\xi(1) - \beta^2\xi(q)\right),$$

where $\xi(q) = q^2/2$. It is easy to check that the infimum is achieved on $q=0$ when $\beta \leq 1$ and $q = 1 - 1/\beta$ when $\beta > 1$. The first case is trivial in many respects, so we will only look at the second case $\beta > 1$ for which the free energy above becomes

$$(1.8) \qquad \mathcal{P}(\beta) = \tfrac{1}{2}(2\beta - \tfrac{3}{2} - \log\beta).$$

We will prove the following bound on $F_N(Q_\varepsilon)$ in (1.6) when the external fields $h_j = 0, j \leq n$. Given a matrix of overlap constraints $Q = (q_{j,j'})$ let us define a matrix

$$(1.9) \qquad \tilde{Q} = (\sqrt{\beta_j\beta_{j'}}q_{j,j'})$$

and let $r_1, \ldots, r_n$ be its eigenvalues. Consider the function

$$(1.10) \qquad f(r) = \begin{cases} \log r + \frac{1}{2}r^2, & \text{for } 0 < r \leq 1, \\ 2r - \frac{3}{2}, & \text{for } 1 \leq r, \end{cases}$$

and note that for $0 < r < 1$,

$$(1.11) \qquad f(r) = \log r + \tfrac{1}{2}r^2 < 2r - \tfrac{3}{2}.$$

The following theorem holds.

THEOREM 1. *For any matrix of overlap constraints $Q$ we have*

$$(1.12) \qquad \limsup_{\varepsilon \to 0} \limsup_{N \to \infty} F_N(Q_\varepsilon) \leq \tfrac{1}{2}\sum_{j \leq n}(f(r_j) - \log\beta_j),$$

*where $(r_j)_{j\leq n}$ are the eigenvalues of (1.9) and $f(r)$ is defined in (1.10). The right-hand side of (1.12) is strictly less than $\sum_{j\leq n}\mathcal{P}(\beta_j)$ if the smallest eigenvalue $\min r_j$ is $<1$.*

The second statement of the theorem follows from (1.11). Indeed, note that if all eigenvalues $r_j$ of $\tilde{Q}$ satisfy $r_j \geq 1$ then the bound (1.12) becomes

$$\sum_{j\leq n} \tfrac{1}{2}(2r_j - \tfrac{3}{2} - \log\beta_j) = \sum_{j\leq n} \tfrac{1}{2}(2\beta_j - \tfrac{3}{2} - \log\beta_j) = \sum_{j\leq n}\mathcal{P}(\beta_j)$$

by (1.8) and since $\sum_{j\leq n} r_j = \text{Tr}(\tilde{Q}) = \sum_{j\leq n}\beta_j$. By (1.11), the bound (1.12) will be strictly less than $\sum_{j\leq n}\mathcal{P}(\beta_j)$ if the smallest eigenvalue $\min r_j < 1$. Let us look at some consequences of this.

EXAMPLE 1. Let $n=2$, $\beta_1, \beta_2 > 1$. For $R_{1,2} \approx u \in [-1,1]$,

(1.13) $$Q = \begin{pmatrix} 1 & u \\ u & 1 \end{pmatrix} \implies \tilde{Q} = \begin{pmatrix} \beta_1 & \sqrt{\beta_1\beta_2}u \\ \sqrt{\beta_1\beta_2}u & \beta_2 \end{pmatrix}$$

and

$$r_1, r_2 = \tfrac{1}{2}(\beta_1 + \beta_2 \pm \sqrt{(\beta_1 - \beta_2)^2 + 4\beta_1\beta_2 u^2}).$$

It is easy to check that $r_2 \geq 1$ if and only if

$$|u| \leq \sqrt{\left(1 - \frac{1}{\beta_1}\right)\left(1 - \frac{1}{\beta_2}\right)} = \sqrt{q_1 q_2}$$

and Theorem 1 (together with standard concentration of measure) implies that the overlap of the coupled system can not exceed $\sqrt{q_1 q_2}$. In the case of two equal temperatures $\beta_1 = \beta_2 = \beta$ this means that the absolute value of the overlap can not exceed $q = 1 - 1/\beta$, which could also be obtained by methods of [8], "breaking" replica symmetric choice of parameters. In Lemma 2 below we will give a general statement that says that if for $j \leq 2$ the overlap of the system $j$ does not exceed $q_j$ than the overlap of a coupled system does not exceed $\sqrt{q_1 q_2}$. Thus, Example 1 could be obtained without the application of Theorem 1. However, (1.12) provides an explicit constructive bound.

Our next example will be less trivial, but before we proceed let us make one observation. Using the fact that for $\beta > 1, h = 0$,

$$\lim_{N \to \infty} F_N(\beta) := \lim_{N \to \infty} \frac{1}{N} \mathbb{E} \log \int_{S_N} \exp \beta H_N(\boldsymbol{\sigma}) \, d\lambda_N(\boldsymbol{\sigma}) = \mathcal{P}(\beta)$$

given by (1.8) and since both $F_N(\beta)$ and $\mathcal{P}(\beta)$ are convex in $\beta$, we have

$$\lim_{N \to \infty} F'_N(\beta) = \lim_{N \to \infty} \frac{1}{2}\beta(1 - \mathbb{E}\langle R_{1,2}^2 \rangle) = \mathcal{P}'(\beta) = 1 - \frac{1}{2\beta},$$

where $\langle \cdot \rangle$ denotes the Gibbs average. This implies that

$$\lim_{N \to \infty} \mathbb{E}\langle R_{1,2}^2 \rangle = \left(1 - \frac{1}{\beta}\right)^2 = q^2.$$

Example 1 for two equal temperatures implies that for any $\varepsilon > 0$,

$$\lim_{N \to \infty} \mathbb{E}\langle I\{|R_{1,2}| \geq q + \varepsilon\}\rangle = 0.$$

These observations combined, of course, imply that for any $\varepsilon > 0$,

(1.14) $$\lim_{N \to \infty} \mathbb{E}\langle I\{|R_{1,2}^2 - q^2| \geq \varepsilon\}\rangle = 0,$$

that is, the overlap can take only values close to $\pm q$.

EXAMPLE 2 (*Ultrametricity*). Let us consider three copies of the system ($n=3$) with the same $\beta > 1$. Ultrametricity means that with high probability in the disorder (randomness of the Hamiltonian $H_N$), for any $\varepsilon > 0$, the Gibbs measure of the event

$$R_{2,3} \geq \min(R_{1,2}, R_{1,3}) - \varepsilon$$

is close to one. By (1.14), the overlaps $R_{j,j'}, j \neq j'$ can only take values $\pm q$ and, thus, the only possible nonultrametric overlap configuration is described by the constraint matrix

$$Q = \begin{pmatrix} 1 & q & q \\ q & 1 & -q \\ q & -q & 1 \end{pmatrix}.$$

It is easy to check that for this matrix

$$r_1 = r_2 = \beta(1+q) \quad \text{and} \quad r_3 = \beta(1-2q).$$

First of all, the matrix $Q$ is positive definite only if $1 - 2q > 0$, that is, $\beta < 2$. Also, for $\beta > 1$ we have $r_3 = \beta(1-2q) = 2 - \beta < 1$ and, therefore, such a configuration $Q$ is not in the support of the Gibbs measure and we have ultrametricity. Some intuition for the fact that the overlap takes two values $\pm q$ is that the Gibbs measure is dominated by two symmetric "states" such that the typical overlap of two spin configurations within each state is equal to $q$ and the overlap of spin configurations from different states is $-q$. Of course, such a picture naturally excludes the overlap configuration given by $Q$ above, but a rigorous proof is another matter.

Let us now go back to the Example 1 for $\beta_1 = \beta_2 = \beta$. For the value of the overlap $R_{1,2} \approx u$ such that $|u| \leq q$, Theorem 1 provides only a trivial bound $2\mathcal{P}(\beta)$ of the type (1.7) while, on the other hand, (1.14) proves that the overlap cannot take values between $-q$ and $q$. At the level of large deviations the bound (1.12) does not detect this and one might ask whether (1.12) is simply not sharp in this case. A similar question may be asked about Example 2 which proves ultrametricity only in the weak sense since due to (1.14) we only had to consider one nonultrametric configuration. However, there are nonultrametric configurations $Q = (q_{j,j'})$, that is, $q_{2,3} < \min(q_{1,2}, q_{1,3})$, for which Theorem 1 will only give a trivial bound $3\mathcal{P}(\beta)$. Again, could it be that (1.12) is simply not sharp in that case? The answer to both of these questions is negative as shown by the following theorem. This result is surprising because it shows that in this model the ultrametricity (and the chaos) cannot be proved at the level of large deviations and, therefore, it is possible that in other models, for example, in the Ising SK model, a similar situation occurs and one should be cautious in one's efforts to prove ultrametricity at the level of free energy.

THEOREM 2. *In the notation, of Theorem 1, if $\beta_j = \beta > 1, j \leq n$ and $\min r_j \geq 1$ then*

$$\lim_{\varepsilon \to 0} \lim_{N \to \infty} F_N(Q_\varepsilon) = n\mathcal{P}(\beta), \tag{1.15}$$

*that is, in this case the bound (1.12) is sharp.*

The proof of this theorem is an extension of the methodology in [7]. It relies on certain a priori estimates, Theorem 6 below, which generally become much more difficult to prove for multiple copies of the system compared to a single system. For example, we do not know how to do this for the Ising SK model or even for the spherical $p$-spin model for even $p > 2$. In the setting of Theorem 2 we were able to prove these estimates using a special "diagonalization" trick developed in Theorem 1.

Finally, we will prove some facts about the overlap of two pure 2-spin systems in the presence of external fields. If $h_j = 0$, we will assume that $\beta_j > 1$ and let $q_j = 1 - 1/\beta_j$. If $h_j \neq 0$, let $\beta_j > 0$ and let $q_j \in [0, 1]$ be a unique solution of

$$h_j^2 + \beta_j^2 q_j = \frac{q_j}{(1-q_j)^2}.$$

Then the following theorem holds.

THEOREM 3. *Let $u_0 = 0$ if $h_1 = 0$ and $h_2 \neq 0$ and let*

$$u_0 = \frac{h_1 h_2 (1-q_1)(1-q_2)}{1 - \beta_1 \beta_2 (1-q_1)(1-q_2)}$$

*if both $h_1 \neq 0, h_2 \neq 0$. Then, for any $\varepsilon > 0$,*

$$\lim_{N \to \infty} \mathbb{E}\langle I\{|R_{1,2} - u_0| \geq \varepsilon\}\rangle = 0.$$

*Pure $p$-spin model, for even $p \geq 4$, without external field.* Next we will consider the pure $p$-spin Hamiltonian (1.2) for even $p \geq 4$, which corresponds to $\xi(q) = q^p/p$, without external field, that is, $h = 0$. It was proven (following an argument of [1]) in Proposition 2.2. in [8] that whenever $\xi''(q)^{-1/2}$ is convex, which is the case here, the infimum in (1.3) or (1.4) is achieved for $k = 2$. When $h = 0$ one can argue that $q_1 = 0$ and, thus, the free energy is

$$\inf_{q,m \in [0,1]} \frac{1}{2} \bigg( \beta^2 \xi(1) + (m-1)\beta^2 \xi(q) \tag{1.16}$$
$$+ \left(1 - \frac{1}{m}\right) \log(1-q) + \frac{1}{m} \log(1 - q(1-m)) \bigg).$$

Proposition 2.3 in [8] states that the infimum will be achieved on $q = 0$ if and only if

(1.17)
$$\sup_{s \leq 0}(\beta^2 \xi(s) + \log(1-s) + s) \leq 0.$$

The case where the infimum is achieved on $q = 0$ corresponds to the trivial case when the overlap can take only the value zero, so, we will only consider the case of $\beta$ large enough, where $q \neq 0$. We will prove the following result concerning the overlap of two pure $p$-spin systems.

THEOREM 4. *Suppose that $\beta_j, j \leq 2$ are such that the infimum in (1.16) is achieved on $q_j \neq 0$, that is, (1.17) fails. Then for any $\varepsilon > 0$,*

(1.18)
$$\lim_{N \to \infty} \mathbb{E}\langle I\{\{|R_{1,2}| \geq \varepsilon\} \cap \{||R_{1,2}| - \sqrt{q_1 q_2}| \geq \varepsilon\}\}\rangle = 0,$$

*that is, the overlap can take only the values $0$ and $\pm\sqrt{q_1 q_2}$.*

The rest of the paper is organized as follows. In Section 2 we prove a general bound on $F_N(Q_\varepsilon)$ using the analogue of Guerra's interpolation. In Section 3 we prove Theorems 1 and 3, in Section 4 we prove Theorem 4 and in Section 5 we prove Theorem 2.

**2. Interpolation.** In this section we will describe the analogue of Guerra's interpolation for the constrained free energy (1.6). Given $k \geq 1$, consider a sequence $\mathbf{m} = (m_l)_{l \leq k}$ such that

$$0 = m_0 < m_1 < \cdots < m_k = 1.$$

We may assume strict inequalities since otherwise in (1.3) and (1.4) we can simply decrease the value of $k$. We consider a sequence for $0 \leq l \leq k+1$ of symmetric $n \times n$ matrices $Q^l = (q^l_{j,j'})_{j,j' \leq n}$ such that $Q^0 = 0$ and such that if we define

(2.1)
$$\hat{Q}^l = (\beta_j \beta_{j'} \xi'(q^l_{j,j'}))$$

then the matrices

(2.2)
$$\Delta_l := \hat{Q}^{l+1} - \hat{Q}^l$$

are nonnegative definite for $0 \leq l \leq k$. Let $\mathbf{z}_l = (z^1_l, z^2_l, \ldots, z^n_l)$ be a Gaussian vector with covariance $\Delta_l$ and let $\mathbf{z}_l$ be independent for $l \leq k$. Finally, let $(\mathbf{z}_{l,i})_{l \leq n}$ be independent copies of $(\mathbf{z}_l)_{l \leq n}$ for $i \leq N$. For $0 \leq t \leq 1$ we define an interpolating Hamiltonian by

(2.3)
$$H_t(\boldsymbol{\sigma}^1, \ldots, \boldsymbol{\sigma}^n) = \sqrt{t} \sum_{j \leq n} \beta_j H_N(\boldsymbol{\sigma}^j)$$
$$+ \sum_{j \leq n} \sum_{i \leq N} \sigma^j_i \left(\sqrt{1-t} \sum_{0 \leq l \leq k} z^j_{l,i} + h_j\right).$$

Let
$$X_{k+1,t}(Q_\varepsilon) = \log \int_{Q_\varepsilon} \exp H_t \, d\lambda_N^n, \tag{2.4}$$

where for simplicity of notations we keep the dependence of $H_t$ and $\lambda_N^n$ on $(\boldsymbol{\sigma}^1, \ldots, \boldsymbol{\sigma}^n)$ implicit. Recursively for $1 \leq l \leq k$ we define

$$X_{l,t}(Q_\varepsilon) = \frac{1}{m_l} \log \mathbb{E}_l \exp m_l X_{l+1,t}(Q_\varepsilon), \tag{2.5}$$

where $\mathbb{E}_l$ denotes the expectation in $(z_{p,i}^j)$ for $p \geq l$. We define

$$\varphi(t) = \frac{1}{N} \mathbb{E} X_{1,t}(Q_\varepsilon). \tag{2.6}$$

Clearly $\varphi(1) = F_N(Q_\varepsilon)$. From now on $\mathcal{R}$ will denote a quantity such that
$$\limsup_{\varepsilon \to 0} \limsup_{N \to \infty} |\mathcal{R}| = 0.$$

The following holds.

THEOREM 5. *We have*
$$\begin{aligned}
\varphi'(t) = {} & \tfrac{1}{2} \sum_{j,j' \leq n} \beta_j \beta_{j'} (\xi(q_{j,j'}) - q_{j,j'} \xi'(q_{j,j'}^{k+1}) + \theta(q_{j,j'}^{k+1})) \\
& - \tfrac{1}{2} \sum_{1 \leq l \leq k} m_l \sum_{j,j' \leq n} \beta_j \beta_{j'} (\theta(q_{j,j'}^{l+1}) - \theta(q_{j,j'}^l)) - R(t) + \mathcal{R},
\end{aligned} \tag{2.7}$$

*where the remainder $R(t) \geq 0$.*

PROOF. The proof of this theorem is a straightforward generalization of Guerra's interpolation for a single system [3] and was explained in detail for coupled copies in [7]. We will not reproduce it here. $\square$

REMARK (*Remainder for $k = 1$*). The remainder $R(t)$ can be written explicitly but we will omit its rather complicated definition in the general case $k \geq 1$ since we will only need the exact form of $R(t)$ in the proof of Theorem 2 for $k = 1$. If $k = 1$, let us define a Hamiltonian
$$h_t(\boldsymbol{\sigma}^1, \ldots, \boldsymbol{\sigma}^n) = \sqrt{t} \sum_{j \leq n} \beta_j H_N(\boldsymbol{\sigma}^j) + \sum_{j \leq n} \sum_{i \leq N} \sigma_i^j (\sqrt{1-t} z_{0,i}^j + h_j) \tag{2.8}$$

and let $\langle \cdot \rangle_t$ denote the Gibbs average on $Q_\varepsilon$ (and its products) with respect to this Hamiltonian. Then given two systems of spins $(\boldsymbol{\sigma}^1, \ldots, \boldsymbol{\sigma}^n) \in Q_\varepsilon$ and $(\boldsymbol{\rho}^1, \ldots, \boldsymbol{\rho}^n) \in Q_\varepsilon$ we have
$$R(t) = \tfrac{1}{2} \sum_{j,j' \leq n} \beta_j \beta_{j'} \mathbb{E} \langle \xi(R^{j,j'}) - R^{j,j'} \xi'(q_{j,j'}^1) + \theta(q_{j,j'}^1) \rangle_t \tag{2.9}$$

where $R^{j,j'} = N^{-1} \sum_{i \leq N} \sigma_i^j \rho_i^{j'}$ is the overlap between $\boldsymbol{\sigma}^j$ and $\boldsymbol{\rho}^{j'}$.

Theorem 5 provides an upper bound for $F_N(Q_\varepsilon)$ since

$$F_N(Q_\varepsilon) = \varphi(1) = \varphi(0) + \int_0^1 \varphi'(t)\,dt$$

and the integral can be bounded using (2.7). Consider a symmetric positive definite $n \times n$ matrix $A = (a_{j,j'})_{j,j' \leq n}$, define $A_{k+1} = A$ and recursively for $l \leq k$ define

(2.10) $$A_l = A_{l+1} - m_l \Delta_l = A_{l+1} - m_l(\hat{Q}^{l+1} - \hat{Q}^l).$$

Below we will always assume that $A$ is chosen so that $A_1$ is also positive definite. We will denote by $|A|$ the determinant of $A$. The following holds.

LEMMA 1. *If* $\mathbf{h} = (h_1, \ldots, h_n)$ *then for any matrix $A$ and the sequence $(A_l)$ defined by (2.10), we have*

(2.11)
$$2\varphi(0) \leq \mathrm{Tr}(A_{k+1}Q) - n + (A_1^{-1}\mathbf{h}, \mathbf{h}) + \mathrm{Tr}(A_1^{-1}\Delta_0)$$
$$+ \sum_{1 \leq l \leq k} \frac{1}{m_l} \log \frac{|A_{l+1}|}{|A_l|} - \log|A_{k+1}| + \mathcal{R}.$$

REMARK. In fact, one can prove that $\lim_{\varepsilon \to 0} \lim_{N \to \infty} 2\varphi(0)$ is equal to the infimum of the right-hand side of (2.11) over all choices of $A$, that is, the bound is sharp. We will give a proof of this statement in the replica symmetric case $k = 1$ which will be needed in the proof of Theorem 2, but the proof of the general case $k \geq 1$ is different only in that the notation is more complicated.

PROOF OF LEMMA 1. To simplify the notation, for $j \leq n$ let

$$\mathbf{z}^j = \left( \sum_{0 \leq l \leq k} z_{l,1}^j + h_j, \ldots, \sum_{0 \leq l \leq k} z_{l,N}^j + h_j \right).$$

Then, since on the set $Q_\varepsilon$ we have $R_{j,j'} = N^{-1}(\boldsymbol{\sigma}^j, \boldsymbol{\sigma}^{j'}) = q_{j,j'} + O(\varepsilon)$, we have

(2.12)
$$\frac{1}{N} X_{k+1,0}(Q_\varepsilon)$$
$$= \frac{1}{N} \log \int_{Q_\varepsilon} \exp \sum_{j \leq n} (\boldsymbol{\sigma}^j, \mathbf{z}^j)\,d\lambda_N^n$$
$$\leq \sum_{j<j'} a_{j,j'} q_{j,j'}$$
$$+ \frac{1}{N} \log \int_{S_N^n} \exp \left( \sum_{j \leq n} (\boldsymbol{\sigma}^j, \mathbf{z}^j) - \sum_{j<j'} a_{j,j'}(\boldsymbol{\sigma}^j, \boldsymbol{\sigma}^{j'}) \right) d\lambda_N^n + O(\varepsilon).$$

We proceed as in Lemma 7.1 in [8]. Let $\nu_{a_{j,j}}(\boldsymbol{\sigma})$ be a Gaussian measure on $\mathbb{R}^N$ of density

$$\left(\frac{a_{j,j}}{2\pi}\right)^{N/2} \exp\left(-\frac{a_{j,j}}{2}\|\boldsymbol{\sigma}\|^2\right).$$

Let us write $\boldsymbol{\rho} \in \mathbb{R}^N$ as $\boldsymbol{\rho} = s\boldsymbol{\sigma}$ where $\boldsymbol{\sigma} \in S_N$ so that by rotational invariance the law of $\boldsymbol{\sigma}$ under $\nu_{a_{j,j}}$ is $\lambda_N$ and $\boldsymbol{\sigma}$ and $s$ are independent. Let $\gamma_j$ be the law of $s$. Let us define $c_j$ by

$$(2.13) \qquad 1 = \gamma_j(\{s \geq c_j\}) \int_{c_j}^{\infty} s \, d\gamma_j(s),$$

which, obviously, implies that $c_j \leq 1$. Therefore,

$$(2.14) \quad \gamma_j(\{s \geq c_j\}) \geq \gamma_j(\{s \geq 1\}) = \nu_{a_{j,j}}(\{\|\boldsymbol{\rho}\| \geq \sqrt{N}\}) = \exp(-N\tau_N(a_{j,j}))$$

and it is easy to check (classical large deviations) that

$$(2.15) \qquad \lim_{N \to \infty} \tau_N(a_{j,j}) = \tau(a_{j,j}) := \tfrac{1}{2}(a_{j,j} - 1 - \log a_{j,j}).$$

By Jensen's inequality and (2.13),

$$\prod_{j \leq n} \gamma_j(\{s_j \geq c_j\}) \int_{S_N^n} \exp\left(\sum_{j \leq n}(\boldsymbol{\sigma}^j, \mathbf{z}^j) - \sum_{j < j'} a_{j,j'}(\boldsymbol{\sigma}^j, \boldsymbol{\sigma}^{j'})\right) d\lambda_N^n$$

$$\leq \int_{\bigcap\{s_j \geq c_j\}} \int_{S_N^n} \exp\left(\sum_{j \leq n}(s_j\boldsymbol{\sigma}^j, \mathbf{z}^j)\right.$$

$$\left. - \sum_{j < j'} a_{j,j'}(s_j\boldsymbol{\sigma}^j, s_{j'}\boldsymbol{\sigma}^{j'})\right) d\lambda_N^n \, d\gamma_1(s_1) \cdots d\gamma_n(s_n)$$

$$(2.16)$$

$$= \int_{\bigcap\{\|\boldsymbol{\rho}^j\| \geq c_j\sqrt{N}\}} \exp\left(\sum_{j \leq n}(\boldsymbol{\rho}^j, \mathbf{z}^j)\right.$$

$$\left. - \sum_{j < j'} a_{j,j'}(\boldsymbol{\rho}^j, \boldsymbol{\rho}^{j'})\right) d\nu_{a_{1,1}}(\boldsymbol{\rho}^1) \cdots d\nu_{a_{n,n}}(\boldsymbol{\rho}^n)$$

$$\leq \int_{\mathbb{R}^{Nn}} \exp\left(\sum_{j \leq n}(\boldsymbol{\rho}^j, \mathbf{z}^j) - \sum_{j < j'} a_{j,j'}(\boldsymbol{\rho}^j, \boldsymbol{\rho}^{j'})\right) d\nu_{a_{1,1}}(\boldsymbol{\rho}^1) \cdots d\nu_{a_{n,n}}(\boldsymbol{\rho}^n).$$

In the last line, all $N$ coordinates are now decoupled. Let $\nu$ be a Gaussian measure on $\mathbb{R}^n$ with covariance $\mathrm{Diag}(a_{1,1}, \ldots, a_{n,n})$ and define for $i \leq N$

$$\mathbf{z}_i = \left(\sum_{l \leq k} z_{l,i}^1 + h_1, \ldots, \sum_{l \leq k} z_{l,i}^n + h_n\right) = \sum_{l \leq k} \mathbf{z}_{l,i} + \mathbf{h}.$$

Then the last line in (2.16) can be rewritten as $\exp \sum_{i \leq N} Y_{k+1,i}$ where

$$
\begin{aligned}
\exp Y_{k+1,i} &= \int_{\mathbb{R}^n} \exp\biggl((\mathbf{x}, \mathbf{z_i}) - \sum_{j<j'} a_{j,j'} x_j x_{j'}\biggr) d\nu(\mathbf{x}) \\
&= \biggl(\frac{\prod a_{j,j}}{(2\pi)^n}\biggr)^{1/2} \int_{\mathbb{R}^n} \exp\biggl((\mathbf{x}, \mathbf{z}_i) - \sum_{j<j'} a_{j,j'} x_j x_{j'} - \frac{1}{2}\sum_{j\leq n} a_{j,j} x_j^2\biggr) d\mathbf{x} \\
&= \biggl(\frac{\prod a_{j,j}}{(2\pi)^n}\biggr)^{1/2} \int_{\mathbb{R}^n} \exp\biggl((\mathbf{x}, \mathbf{z_i}) - \frac{1}{2}(A\mathbf{x}, \mathbf{x})\biggr) d\mathbf{x} \\
&= \biggl(\frac{\prod a_{j,j}}{|A|}\biggr)^{1/2} \exp\biggl(\frac{1}{2}(A^{-1}\mathbf{z}_i, \mathbf{z}_i)\biggr).
\end{aligned}
$$

(2.17)

Combining this with (2.12), (2.14), (2.15) and (2.16) we get

$$
\begin{aligned}
\frac{1}{N} & X_{k+1,0}(Q_\varepsilon) \\
&\leq \sum_{j<j'} a_{j,j'} q_{j,j'} + \sum_{j\leq n} \tau(a_{j,j}) + \frac{1}{2}\log\frac{\prod a_{j,j}}{|A|} \\
&\quad + \frac{1}{N}\sum_{i\leq N} \frac{1}{2}(A^{-1}\mathbf{z}_i, \mathbf{z}_i) + \mathcal{R} \\
&= \frac{1}{2}(\mathrm{Tr}(AQ) - n - \log|A|) + \frac{1}{N}\sum_{i\leq N} \frac{1}{2}(A^{-1}\mathbf{z}_i, \mathbf{z}_i) + \mathcal{R},
\end{aligned}
$$

(2.18)

where $\mathcal{R} = O(\varepsilon) + c_N$ for some $c_N \to 0$. This bound will propagate in the recursion (2.5) and since $\mathbf{z}_i$ are all independent we only need to compute what happens with one of the terms $(A^{-1}\mathbf{z}_i, \mathbf{z}_i)/2$. Namely, if

$$
y_{k+1} = \frac{1}{2}(A^{-1}\mathbf{z}_i, \mathbf{z}_i), \qquad y_l = \frac{1}{m_l}\log \mathbb{E}_l \exp m_l y_{l+1} \qquad \text{for } l \leq k,
$$

then by induction over $l$ it should be obvious that

(2.19) $\qquad \varphi(0) \leq \tfrac{1}{2}(\mathrm{Tr}(AQ) - n - \log|A|) + \mathbb{E} y_1 + \mathcal{R}.$

To compute the sequence $(y_l)$ we use the fact that for a Gaussian vector $\mathbf{g}$ with covariance $C$, parameter $m > 0$ and a symmetric positive definite matrix $A$ such that $A - mC$ is also positive definite, we have

$$
\begin{aligned}
\frac{1}{m}\log \mathbb{E}\exp \frac{m}{2}&(A^{-1}(\mathbf{x}+\mathbf{g}), \mathbf{x}+\mathbf{g}) \\
&= \frac{1}{2m}\log\frac{|A|}{|A-mC|} + \frac{1}{2}((A-mC)^{-1}\mathbf{x}, \mathbf{x}),
\end{aligned}
$$

which is a simple exercise. Since the covariance of $\mathbf{z}_{l,i}$ is $\Delta_l$, using (2.10) we get by induction over $l \leq k$ that

$$y_1 = \frac{1}{2}(A_1^{-1}(\mathbf{z}_{0,i} + \mathbf{h}), \mathbf{z}_{0,i} + \mathbf{h}) + \frac{1}{2} \sum_{1 \leq l \leq k} \frac{1}{m_l} \log \frac{|A_{l+1}|}{|A_l|}$$

and, therefore,

$$\mathbb{E}y_1 = \frac{1}{2}\operatorname{Tr}(A_1^{-1}\Delta_0) + \frac{1}{2}(A_1^{-1}\mathbf{h}, \mathbf{h}) + \frac{1}{2} \sum_{1 \leq l \leq k} \frac{1}{m_l} \log \frac{|A_{l+1}|}{|A_l|}.$$

Plugging this back into (2.19) completes the proof. □

Combining Theorem 5 and Lemma 1 gives an upper bound on $F_N(Q_\varepsilon)$. From now on we will always set

(2.20) $$Q^{k+1} = Q$$

in which case the first sum on the right-hand side of (2.7) disappears and the bound becomes

(2.21) $$\begin{aligned} 2F_N(Q_\varepsilon) \leq {} & \operatorname{Tr}(A_{k+1}Q) - n + (A_1^{-1}\mathbf{h}, \mathbf{h}) + \operatorname{Tr}(A_1^{-1}\Delta_0) \\ & + \sum_{1 \leq l \leq k} \frac{1}{m_l} \log \frac{|A_{l+1}|}{|A_l|} - \log |A_{k+1}| \\ & - \sum_{1 \leq l \leq k} m_l \sum_{j,j' \leq n} \beta_j \beta_{j'} (\theta(q_{j,j'}^{l+1}) - \theta(q_{j,j'}^l)) + \mathcal{R}. \end{aligned}$$

It will often be convenient to think of $A_k$ as a free parameter and think of $A_{k+1}$ as $A_k + \Delta_k$ and write (2.21) as

(2.22) $$\begin{aligned} 2F_N(Q_\varepsilon) \leq {} & \operatorname{Tr}(\Delta_k Q) + \operatorname{Tr}(A_k Q) - n + (A_1^{-1}\mathbf{h}, \mathbf{h}) \\ & + \operatorname{Tr}(A_1^{-1}\Delta_0) + \sum_{1 \leq l \leq k-1} \frac{1}{m_l} \log \frac{|A_{l+1}|}{|A_l|} - \log |A_k| \\ & - \sum_{1 \leq l \leq k} m_l \sum_{j,j' \leq n} \beta_j \beta_{j'} (\theta(q_{j,j'}^{l+1}) - \theta(q_{j,j'}^l)) + \mathcal{R}. \end{aligned}$$

In general, it is not clear how to choose the parameters in this bound that would witness the trivial bound (1.7). In the next section we will show how this can be done for the pure 2-spin model.

## 3. Pure 2-spin model.

PROOF OF THEOREM 1. We will first prove Theorem 1 by constructing the sequence $(Q^l)$ and $A$ in a special way. Let us mention that in order to prove Theorem 1 it would be enough to consider the bound (2.21) only for $k = 1$. However, to illustrate the general idea we will look at any $k \geq 1$ and at the end of our argument it will be clear why for 2-spin model one should take $k = 1$. Note that now $\xi(q) = \theta(q) = q^2/2$ and $\xi'(q) = q$. Let

$$\tilde{Q}^{k+1} = (\tilde{q}^{k+1}_{j,j'}) := (\sqrt{\beta_j \beta_{j'}} q^{k+1}_{j,j'}) = (\sqrt{\beta_j \beta_{j'}} q_{j,j'})$$

and let

$$\tilde{Q}^{k+1} = O^T R^{k+1} O$$

be its Jordan decomposition for some $R^{k+1} = \text{Diag}(r_1, \ldots, r_n)$ and orthogonal matrix $O$. We now consider a nondecreasing sequence

$$0 = R^0 \leq R^1 \leq \cdots \leq R^k \leq R^{k+1}$$

of diagonal matrices $R^l = \text{Diag}(r^l_1, \ldots, r^l_n)$ which means that each sequence $(r^l_j)_{l \leq k+1}$ is nondecreasing and define

(3.1) $$\tilde{Q}^l = (\tilde{q}^l_{j,j'}) = (\sqrt{\beta_j \beta_{j'}} q^l_{j,j'}) = O^T R^l O,$$

which also defines $Q^l$ and $\hat{Q}^l$, in particular,

$$\hat{q}^l_{j,j'} = \beta_j \beta_{j'} q^l_{j,j'} = \sqrt{\beta_j \beta_{j'}} \tilde{q}^l_{j,j'}.$$

Next, for $B = \text{Diag}(b_1, \ldots, b_n)$ where each $b_j > 0$ we let

(3.2) $$D = (d_{j,j'}) = O^T B O \quad \text{and} \quad A = (\sqrt{\beta_j \beta_{j'}} d_{j,j'}).$$

Then by (2.10) the elements of the matrix $A_l$ are

$$\sqrt{\beta_j \beta_{j'}} \left( d_{j,j'} - \sum_{l \leq p \leq k} m_p(\tilde{q}^{p+1}_{j,j'} - \tilde{q}^p_{j,j'}) \right)$$

$$= \sqrt{\beta_j \beta_{j'}} \left( D - \sum_{l \leq p \leq k} m_p(\tilde{Q}^{p+1} - \tilde{Q}^p) \right)_{j,j'}$$

$$= \sqrt{\beta_j \beta_{j'}} \left( O^T \left( B - \sum_{l \leq p \leq k} m_p(R^{p+1} - R^p) \right) O \right)_{j,j'}.$$

It is obvious that for any matrix $M = (m_{j,j'})$ we have

$$|(\sqrt{\beta_j \beta_{j'}} m_{j,j'})| = |M| \prod_{j \leq n} \beta_j$$

and, therefore,

(3.3) $$|A_l| = |B_l| \prod_{j \leq n} \beta_j \quad \text{for } B_l = B - \sum_{l \leq p \leq k} m_p(R^{p+1} - R^p).$$

Similarly,

$$\text{Tr}(A_1^{-1}\Delta_0) = \lim_{m \to 0} \frac{1}{m} \log \frac{|A_1|}{|A_1 - m\Delta_0|}$$

(3.4)

$$= \lim_{m \to 0} \frac{1}{m} \log \frac{|B_1|}{|B_1 - mR^1|} = \text{Tr}(B_1^{-1}R^1) = \sum_{j \leq n} \frac{r_j^1}{b_j^1},$$

where we denoted by $b_j^1$ the $j$th element on the diagonal of $B_1$. Next,

$$\text{Tr}(A_{k+1}Q) = \sum_{j,j'} a_{j,j'} q_{j,j'} = \sum_{j,j'} \sqrt{\beta_j \beta_{j'}} d_{j,j'} q_{j,j'} = \text{Tr}(D\tilde{Q}^{k+1})$$

(3.5)

$$= \text{Tr}(O^T B_{k+1} O O^T R^{k+1} O) = \text{Tr}(BR^{k+1}) = \sum_{j \leq n} b_j r_j.$$

Finally,

$$\sum_{l \leq k} m_l \sum_{j,j' \leq n} \beta_j \beta_{j'} (\theta(q_{j,j'}^{l+1}) - \theta(q_{j,j'}^l))$$

(3.6)

$$= \tfrac{1}{2} \sum_{l \leq k} m_l \sum_{j,j' \leq n} ((\tilde{q}_{j,j'}^{l+1})^2 - (\tilde{q}_{j,j'}^l)^2)$$

$$= \tfrac{1}{2} \sum_{l \leq k} m_l \sum_{j \leq n} ((r_j^{l+1})^2 - (r_j^l)^2) = \sum_{l \leq k} m_l \sum_{j \leq n} (\theta(r_j^{l+1}) - \theta(r_j^l)),$$

since $\tilde{Q}^l = O^T R^l O$ and the Frobenius norm of the matrix is preserved by an orthogonal transformation, that is,

$$\sum_{j,j' \leq n} (\tilde{q}_{j,j'}^l)^2 = \sum_{j \leq n} (r_j^l)^2.$$

When $\mathbf{h} = 0$, (3.3)–(3.6) result in a bound (2.21) being transformed into

$$2F_N(Q_\varepsilon) \leq \text{Tr}(BR^{k+1}) - n + \text{Tr}(B_1^{-1}R^1) + \sum_{1 \leq l \leq k} \frac{1}{m_l} \log \frac{|B_{l+1}|}{|B_l|}$$

(3.7)

$$- \log |B| - \log \prod_{j \leq n} \beta_j - \sum_{1 \leq l \leq k} m_l \sum_{j \leq n} (\theta(r_j^{l+1}) - \theta(r_j^l)) + \mathcal{R}.$$

Obviously, if we denote the diagonal elements of $B_l$ by $b_1^l, \ldots, b_n^l$ so that

$$b_j^l = b_j - \sum_{l \leq p \leq k} m_p(r_j^{p+1} - r_j^p)$$

then (3.7) decouples into the sum of $n$ terms

$$\frac{r_j^1}{b_j^1} + b_j r_j^{k+1} - 1 - \log b_j - \log \beta_j + \sum_{1 \leq l \leq k} \frac{1}{m_l} \log \frac{b_j^{l+1}}{b_j^l}$$

(3.8)
$$-\sum_{1\leq l\leq k} m_l(\theta(r_j^{l+1}) - \theta(r_j^l))$$

which was our main idea and motivation. Note that the same idea does not immediately work for $p$-spin model for $p > 2$ since (3.6) fails there. We will now explain why minimizing this bound it is enough to look only at the case $k = 1$. Besides the explicit term $-\log \beta_j$, (3.8) is very similar to (1.3) for $\beta = 1$ with minor differences, namely, the range of the sequence $(r_j^l)_{l\leq k}$ is between 0 and $r_j^{k+1}$—the $j$th eigenvalue of the matrix $\tilde{Q}^{k+1}$—and we have $b_j r_j^{k+1}$ instead of $b$. Similar to (1.4), one can repeat the proof of the Crisanti–Sommers formula in Section 4 in [8] to show that the infimum of (3.8) over $b_j$ and the sequence $(r_j^l)_{l\leq k}$ is equal to the infimum of

(3.9) $\quad \dfrac{r_j^1}{\delta_j^1} + \sum_{1\leq l\leq k-1} \dfrac{1}{m_l}\log\dfrac{\delta_j^l}{\delta_j^{l+1}} + \log\delta_j^k + \sum_{l\leq k} m_l(\xi(r_j^{l+1}) - \xi(r_j^l)) - \log\beta_j,$

where $\delta_j^l = \sum_{l\leq p\leq k} m_p(r_j^{p+1} - r_j^p)$. (This is especially easy to check for $k = 1$.) The proof of Proposition 2.2 in [8] yields that the infimum of (3.9) is achieved on the replica symmetric choice of parameters when $k = 1$,

$$\inf_{r_j^1 \in [0, r_j^2]} \left( \log(r_j^2 - r_j^1) + \dfrac{r_j^1}{r_j^2 - r_j^1} + \dfrac{1}{2}((r_j^2)^2 - (r_j^1)^2) - \log\beta_j \right).$$

It is easy to check that the infimum is achieved on

(3.10) $\qquad\qquad r_j^1 = \begin{cases} r_j^2 - 1, & \text{if } r_j^2 \geq 1, \\ 0, & \text{if } r_j^2 < 1, \end{cases}$

and the infimum is equal to $f(r_j^2) - \log\beta_j$ which proves Theorem 1. $\square$

REMARK. In Theorem 2 we will prove that the bound of Theorem 1 is sharp in when $\beta_j = \beta > 1$ for $j \leq n$ and when the bound is $n\mathcal{P}(\beta)$. The proof of Theorem 1 shows that this happens when all $r_j^2 \geq 1$ and (3.10) implies that $r_j^1 = r_j^2 - 1$. Definition (3.1) then gives that

$$Q^2 - Q^1 = \beta^{-1} O^T (R^2 - R^1) O = \beta^{-1} I$$

and

(3.11) $\qquad Q^1 = Q^2 - \beta^{-1}I \quad \text{and} \quad \Delta_1 = \beta^2(Q^2 - Q^1) = \beta I.$

The diagonal elements of $Q^1$ are $q_{j,j}^1 = 1 - \beta^{-1} = q$ and off-diagonal element are $q_{j,j'}^1 = q_{j,j'}^2 = q_{j,j'}$. It is also easy to check that the infimum in (3.8) is achieved for $b_j = 2$, that is, $B = 2I$ and definition (3.2) implies that $A = 2\beta I$. In fact, instead of checking that the proof of Theorem 1 produces these parameters $Q^1$ and $A$, one could simply use these parameters in (2.21) and see that they result in a bound $F_N(Q_\varepsilon) \leq n\mathcal{P}(\beta) + \mathcal{R}$. Clearly, the condition $\min r_j \geq 1$ that gives the bound $n\mathcal{P}(\beta)$ is equivalent to saying that $Q^1 = Q - \beta^{-1}I$ is nonnegative definite.

PROOF OF THEOREM 3. We already mentioned that for $\xi(q) = q^2/2$ the infimum in (1.4) [and (1.3)] is achieved for $k = 1$,

$$\mathcal{P}(\beta, h) = \inf_{q \in [0,1]} \frac{1}{2}\Big(h^2(1-q) + \frac{q}{1-q} \tag{3.12}$$
$$+ \log(1-q) + \beta^2\xi(1) - \beta^2\xi(q)\Big).$$

The critical point condition for $q$ is

$$h^2 + \beta^2\xi'(q) = h^2 + \beta^2 q = \frac{q}{(1-q)^2}. \tag{3.13}$$

If $h_j \neq 0$, let $q_j$ be the unique solution of (3.13) that corresponds to the choice of parameters $\beta_j, h_j$. If $h_j = 0$, we will assume as before that $\beta_j > 1$ and $q_j = 1 - 1/\beta_j$. For $k = 1$, (2.22) is

$$2F_N(Q_\varepsilon) \leq \mathrm{Tr}(A_1 Q) + \mathrm{Tr}(\Delta_1 Q) - 2 + (A_1^{-1}\mathbf{h}, \mathbf{h}) + \mathrm{Tr}(A_1^{-1}\Delta_0)$$
$$- \log|A_1| - \sum_{j,j' \leq 2} \beta_j \beta_{j'}(\theta(q_{j,j'}) - \theta(q^1_{j,j'})) + \mathcal{R}.$$

When $R_{1,2} \approx q_{1,2} = u$ for $|u| \leq \sqrt{q_1 q_2}$, we will take

$$Q^1 = \begin{pmatrix} q_1 & u \\ u & q_2 \end{pmatrix} \quad \text{and} \quad A_1 = \begin{pmatrix} a_1 & \lambda \\ \lambda & a_2 \end{pmatrix},$$

and the bound becomes $2F_N(Q_\varepsilon) \leq U + \mathcal{R}$ for

$$U := \frac{\beta_1^2}{2}(1-q_1)^2 + \frac{\beta_2^2}{2}(1-q_2)^2 + a_1 + a_2 - 2 + 2\lambda u - \log(a_1 a_2 - \lambda^2) \tag{3.14}$$
$$+ \frac{1}{a_1 a_2 - \lambda^2}(a_2(\beta_1^2 q_1 + h_1^2) + a_1(\beta_2^2 q_2 + h_2^2) - 2(\beta_1\beta_2 u + h_1 h_2)\lambda).$$

It is easy to see that for $a_j = 1/(1-q_j) = \beta_j$ and $\lambda = 0$

$$\tfrac{1}{2}U = \mathcal{P}(\beta_1, h_1) + \mathcal{P}(\beta_2, h_2). \tag{3.15}$$

The derivative of (3.14) in $\lambda$ at $\lambda = 0$ and $a_j = 1/(1-q_j) = \beta_j$ is

$$\frac{1}{2}\frac{\partial U}{\partial \lambda}\bigg|_{\lambda=0} = u - \frac{1}{a_1 a_2}(\beta_1\beta_2 u + h_1 h_2) = u - (1-q_1)(1-q_2)(\beta_1\beta_2 u + h_1 h_2).$$

If both $h_j = 0$ then since $\beta_j = 1/(1-q_j)$ this derivative is equal to 0 for all $|u| \leq \sqrt{q_1 q_2}$ which means that we can not improve upon (3.15) by fluctuating $\lambda$ around zero. Suppose that $h_1 = 0$ and $h_2 \neq 0$. Then again $\beta_1 = 1/(1 - q_1)$ and by (3.13), $\beta_2 \neq 1/(1 - q_2)$. The derivative

$$\frac{1}{2}\frac{\partial U}{\partial \lambda}\bigg|_{\lambda=0} = u(1 - \beta_2(1 - q_2))$$

is equal to zero only if $u = u_0 = 0$. Therefore, for $|u| \leq \sqrt{q_1 q_2}$ and $u \neq 0$, by fluctuating $\lambda$ around zero we can prove that $F_N(Q_\varepsilon) < \mathcal{P}(\beta_1, h_1) + \mathcal{P}(\beta_2, h_2) + \mathcal{R}$ and by concentration of measure this means that the overlap can take only the value 0 between $\pm\sqrt{q_1 q_2}$. If both $h_j \neq 0$ then (3.13) implies that $\beta_j(1 - q_j) < 1$ and thus the derivative is equal to zero only at one point

$$u_0 = \frac{h_1 h_2 (1 - q_1)(1 - q_2)}{1 - \beta_1 \beta_2 (1 - q_1)(1 - q_2)}$$

and again the overlap can take only the value $u_0$ between $\pm\sqrt{q_1 q_2}$. To finish the proof of Theorem 3 we need to show that the overlap can not take values $|u| > \sqrt{q_1 q_2}$ which will be done in Lemma 2 below. $\square$

LEMMA 2. *For $j \leq 2$ let $\langle \cdot \rangle_j$ denote the Gibbs average of the (random) Gibbs measure $G_j$ on $S_N$ (and its products) and let $\langle \cdot \rangle_{1,2}$ denote the Gibbs average with respect to the product Gibbs measure $G_1 \otimes G_2$ on $S_N^2$. If for any $\varepsilon > 0$*

$$\lim_{N \to \infty} \mathbb{E}\langle I\{|R_{1,2}| \geq q_j + \varepsilon\}\rangle_j = 0$$

*then for any $\varepsilon > 0$*

$$\lim_{N \to \infty} \mathbb{E}\langle I\{|R_{1,2}| \geq \sqrt{q_1 q_2} + \varepsilon\}\rangle_{1,2} = 0.$$

PROOF. For any integer $k \geq 1$ we can write

$$\mathbb{E}\langle R_{1,2}^k \rangle_{1,2} = N^{-k} \sum_{i_1,\ldots,i_k} \mathbb{E}\langle \sigma_{i_1}^1 \cdots \sigma_{i_k}^1 \sigma_{i_1}^2 \cdots \sigma_{i_k}^2 \rangle_{1,2}$$

$$= N^{-k} \sum_{i_1,\ldots,i_k} \mathbb{E}\langle \sigma_{i_1}^1 \cdots \sigma_{i_k}^1 \rangle_1 \langle \sigma_{i_1}^2 \cdots \sigma_{i_k}^2 \rangle_2$$

$$\leq N^{-k} \sum_{i_1,\ldots,i_k} (\mathbb{E}\langle \sigma_{i_1}^1 \cdots \sigma_{i_k}^1 \rangle_1^2)^{1/2}(\mathbb{E}\langle \sigma_{i_1}^2 \cdots \sigma_{i_k}^2 \rangle_2^2)^{1/2}$$

$$\leq \left(N^{-k} \sum_{i_1,\ldots,i_k} \mathbb{E}\langle \sigma_{i_1}^1 \cdots \sigma_{i_k}^1 \rangle_1^2\right)^{1/2} \left(N^{-k} \sum_{i_1,\ldots,i_k} \mathbb{E}\langle \sigma_{i_1}^2 \cdots \sigma_{i_k}^2 \rangle_2^2\right)^{1/2}$$

$$= \sqrt{\mathbb{E}\langle R_{1,2}^k \rangle_1 \mathbb{E}\langle R_{1,2}^k \rangle_2}$$

and the result is obvious. $\square$

REMARK. Another approach to proving the above results, as well as Theorem 2, for $p = 2$ model would be diagonalizing the system from the beginning by writing the Hamiltonian after orthogonal transformation as $\sum_{i \leq N} \lambda_i \sigma_i^2$, where $\lambda_i$ are the eigenvalues of the matrix $(g_{ij})_{i,j \leq N}$. In [2], a large deviation principle was proved for

$$\left( \frac{1}{N} \sum_{i \leq N} \sigma_i^2, \frac{1}{N} \sum_{i \leq N} \lambda_i \sigma_i^2 \right)$$

and one can similarly prove large deviation principles for several copies of the system by including a vector of overlaps. This approach would, probably, yield simpler proofs of our results for $p = 2$ model but in the present paper we are trying to develop a methodology that could potentially be used for $p > 2$ models or mixed $p$-spin models in which case the diagonalization idea would not work. In the next section we apply our methodology to the case of two copies, $n = 2$, of the pure $p$-spin model for $p > 2$.

**4. Pure $p$-spin model.** In this section we will prove Theorem 4. Unfortunately, the argument of the previous section does not work because for $p > 2$ equation (3.6) does not hold anymore. Consequently, we were unable to prove a general result for all $n \geq 2$. However, for pure $p$-spin model without external field the parameters $(m, q)$ that achieve a minimum in the free energy formula (1.16) satisfy certain special properties that will give us enough information to find the informative parameters in the bound (2.22) for $n = 2$. These properties are given in Lemma 3 below.

We assume that $h = 0$ and in order to avoid trivial situations we consider only $\beta > 0$ such that (1.17) fails. First we will gather important information about the parameters $(m, q)$ that achieve the infimum in (1.16). The critical point conditions for $(m, q)$ are given by

$$\beta^2 \xi'(q) = \frac{1}{m} \left( \frac{1}{1-q} - \frac{1}{1-q+mq} \right),$$

(4.1)

$$\beta^2 \xi(q) = \frac{1}{m^2} \log\left( \frac{1-q+qm}{1-q} \right) - \frac{q}{m} \frac{1}{1-q+qm}$$

and from now on we assume that $(m, q)$ satisfy (4.1). Let $x > 0$ be a unique solution of

(4.2) $$\frac{1}{p} = \frac{1+x}{x^2} \log(1+x) - \frac{1}{x}.$$

It is easy to check that the right-hand side is a convex function decreasing from $1/2$ to $0$ as $x$ increases from $0$ to $\infty$ and, hence, for $p \geq 2$ there is, indeed, a unique solution $x$. Define

(4.3) $$\delta = \frac{x}{1+x}, \qquad \gamma = \frac{p-1}{p} \frac{x^2}{1+x}.$$

Then the following holds.

LEMMA 3. *For all $\beta$ such that (1.17) fails we have*

(4.4) $\quad qm = x(1-q), \qquad m\beta^2 \xi'(q)(1-q) = \delta, \qquad m^2\beta^2\theta(q) = \gamma.$

PROOF. Introducing the notation $A = 1-q, B = qm$, (4.1) can be rewritten as

$$\beta^2 \frac{\xi'(q)}{q} = \frac{1}{B}\left(\frac{1}{A} - \frac{1}{A+B}\right) = \frac{1}{A(A+B)},$$

$$\beta^2 \frac{\xi(q)}{q^2} = \frac{1}{B^2}\log\left(1 + \frac{B}{A}\right) - \frac{1}{B(A+B)}.$$

Since for $\xi(q) = q^p/p$ we have

$$\frac{\xi(q)}{q^2} = \frac{1}{p}\frac{\xi'(q)}{q},$$

the above equations imply that

$$\frac{1}{pA(A+B)} = \frac{1}{B^2}\log\left(1 + \frac{B}{A}\right) - \frac{1}{B(A+B)}.$$

This shows that $B/A$ satisfies (4.2) and, therefore, $B = xA$. From the equations above we get

(4.5) $\qquad m\beta^2 \xi'(q)(1-q) = \dfrac{B}{A+B} = \dfrac{x}{1+x} = \delta$

and since $\theta(q) = (1 - 1/p)q\xi'(q)$,

(4.6)
$$m^2\beta^2\theta(q) = \frac{p-1}{p}m^2 q\beta^2 \xi'(q) = \frac{p-1}{p}B\left(\frac{1}{A} - \frac{1}{A+B}\right)$$
$$= \frac{p-1}{p}\frac{x^2}{1+x} = \gamma. \qquad \square$$

Even though we established Lemma 2 by looking at the Crisanti–Sommers representation (1.4), we will use (1.3) to write $\mathcal{P}(\beta) := \mathcal{P}(\beta, 0)$ in order to compare more easily with the bound (2.21), which is not written in the Crisanti–Sommers form. [Of course, it would be possible to write (2.21) in the Crisanti–Sommers form with extra work.] For $k = 2$ and $q_1 = 0, q_2 = q$, (1.3) becomes

(4.7)
$$2\mathcal{P}(\beta) = \inf_{b,m,q}\Big(b - 1 + \frac{1}{m}\log\frac{D_2}{D_1}$$
$$- \log D_2 - m\beta^2\theta(q) - \beta^2(\theta(1) - \theta(q))\Big).$$

Since $D_2 = b - \beta^2(\xi'(1) - \xi'(q))$ we can minimize over $D_2$ instead of $b$ and to simplify the notations we will simply make the change of variable
$$b \to b + \beta^2(\xi'(1) - \xi'(q)).$$
Then, since $D_1 = D_2 - m\beta^2\xi'(q)$, we have

$$2\mathcal{P}(\beta) = \inf_{b,m,q} \Big(\beta^2\xi'(1) - \beta^2\xi'(q) + b - 1 - \log b \tag{4.8}$$
$$- \frac{1}{m}\log\Big(1 - \frac{1}{b}m\beta^2\xi'(q)\Big) - m\beta^2\theta(q) - \beta^2(\theta(1) - \theta(q))\Big).$$

The infimum is achieved on $(m, q)$ as in (4.1) and $b = 1/(1-q)$. Using (4.4), we get

$$2\mathcal{P}(\beta) = \beta^2\xi'(1) - \beta^2\xi'(q) + b - 1 - \log b - \frac{1}{m}\log(1-\delta) - \frac{\gamma}{m} \tag{4.9}$$
$$- \beta^2(\theta(1) - \theta(q)).$$

PROOF OF THEOREM 4. Given $\beta_j, j \leq 2$ as in Theorem 4, let $m_j$ and $q_j$ be the corresponding solutions of (4.1) and $b_j = 1/(1-q_j)$ so that $\mathcal{P}(\beta_j)$ is given by (4.9). Let us first consider the overlap $R_{1,2} \approx u$ for $|u| \leq \sqrt{q_1 q_2}$. Define $c$ by $|u| = c\sqrt{q_1 q_2}$. We are going to use the bound (2.22) for $k=4$ with $\mathbf{m} = (0, m, m_1, m_2, 1)$ for some $m > 0$ that will be specified later and

$$A_4 = \begin{pmatrix} b_1 & 0 \\ 0 & b_2 \end{pmatrix}, \quad Q^1 = \begin{pmatrix} cq_1 & u \\ u & cq_2 \end{pmatrix},$$
$$Q^2 = \begin{pmatrix} q_1 & u \\ u & cq_2 \end{pmatrix}, \quad Q^3 = \begin{pmatrix} q_1 & u \\ u & q_2 \end{pmatrix}.$$

With these parameters the right-hand side of (2.22) becomes $U(m,c) + \mathcal{R}$ for

$$U(m,c) := -\beta_1^2(\theta(1) - \theta(q_1)) - \beta_2^2(\theta(1) - \theta(q_2)) + \beta_1^2(\xi'(1) - \xi'(q_1))$$
$$+ \beta_2^2(\xi'(1) - \xi'(q_2))$$
$$+ b_1 - 1 - \log b_1 + b_2 - 1 - \log b_2 + \mathrm{I} + \mathrm{II},$$

where

$$\mathrm{I} = -(m_1\beta_1^2\theta(q_1) + m_2\beta_2^2\theta(q_2))(1 - c^p)$$
$$- m(\beta_1^2\theta(q_1) + \beta_2^2\theta(q_2) + 2\beta_1\beta_2\theta(\sqrt{q_1 q_2}))c^p$$

and

$$\mathrm{II} = -\frac{1}{m_1}\log\Big(1 - \frac{1}{b_1}m_1\beta_1^2\xi'(q_1)(1 - c^{p-1})\Big)$$
$$- \frac{1}{m_2}\log\Big(1 - \frac{1}{b_2}m_2\beta_2^2\xi'(q_2)(1 - c^{p-1})\Big) - \frac{1}{m}\log(1 - mS)$$

for
$$S = \sum_{j \leq 2} \frac{\beta_j^2 \xi'(q_j) c^{p-1}}{b_j - m_j \beta_j^2 \xi'(q_j)(1 - c^{p-1})}.$$

Since $\theta(\sqrt{q_1 q_2}) = \sqrt{\theta(q_1)\theta(q_2)}$, Lemma 4.4 implies that
$$\mathrm{I} = -\left(\frac{1}{m_1} + \frac{1}{m_2}\right)\gamma(1 - c^p) - m\left(\frac{1}{m_1} + \frac{1}{m_2}\right)^2 \gamma c^p$$

and
$$\mathrm{II} = -\frac{1}{m_1}\log(1 - \delta(1 - c^{p-1})) - \frac{1}{m_2}\log(1 - \delta(1 - c^{p-1}))$$
$$- \frac{1}{m}\log\left(1 - m\left(\frac{1}{m_1} + \frac{1}{m_2}\right)\frac{\delta c^{p-1}}{1 - \delta(1 - c^{p-1})}\right).$$

If we define $m_0$ by

(4.10) $$\frac{1}{m_0} = \frac{1}{m_1} + \frac{1}{m_2}$$

then for $m = m_0$
$$\mathrm{I} + \mathrm{II} = -\left(\frac{1}{m_1} + \frac{1}{m_2}\right)\gamma - \left(\frac{1}{m_1} + \frac{1}{m_2}\right)\log(1 - \delta).$$

Therefore, comparing with (4.9), for all $c$
$$U(m_0, c) = 2\mathcal{P}(\beta_1) + 2\mathcal{P}(\beta_2).$$

The derivative
$$d(c) = \left.\frac{\partial U(m, c)}{\partial m}\right|_{m=m_0} = \frac{1}{m_0^2}\left(-\gamma c^p + \log\frac{1 - \delta}{1 - \delta(1 - c^{p-1})} + \frac{\delta c^{p-1}}{1 - \delta}\right)$$
$$= \frac{1}{m_0^2}\left(-\frac{p-1}{p}\frac{x^2}{1+x}c^p - \log(1 + xc^{p-1}) + xc^{p-1}\right)$$

satisfies

(4.11)    $d(0) = d(1) = 0$ and $d(c) < 0$    for $0 < c < 1$.

The equality $d(1) = 0$ is equivalent to (4.2) and the fact that $d(c) < 0$ for $0 < c < 1$ is a consequence of the following. The derivative $d'(c)$ is equal to 0 if $c = 0, c = 1$ or if
$$1 + xc^{p-1} = (1 + x)c^{p-2}.$$

Making a change of variable $a = c^{p-2}$ gives
$$1 + xa^{(p-1)/(p-2)} = (1 + x)a.$$

The left-hand side is convex in $a$, the right-hand side is linear and two sides are equal at $a = 1$. Therefore, there is at most one solution on the interval $(0,1)$ and $d'(c) = 0$ has at most one solution on $(0,1)$. However, since $d(0) = d(1) = 0$, $d'(c) = 0$ has exactly one solution on $(0,1)$. Thus, $d(c)$ cannot change sign inside $(0,1)$ which proves (4.11). For $0 < |u| < \sqrt{q_1 q_2}$, fluctuating $m$ near $m_0$ proves that

$$F_N(Q_\varepsilon) \leq \tfrac{1}{2} \inf_m U(m, c) + \mathcal{R} < \mathcal{P}(\beta_1) + \mathcal{P}(\beta_2) + \mathcal{R}.$$

Therefore, by concentration of measure the overlap cannot take values $u$ with $0 < |u| < \sqrt{q_1 q_2}$. Lemma 2 will finish the proof of Theorem 4 if we can show that the overlap of two systems at the same $\beta$ cannot take values $|u| > q = 1 - 1/\beta$. This can be shown by the usual replica symmetry breaking as in [8], which we will explain here for completeness. For specificity, let us assume that $u > q$. We are going to use the bound (2.22) for $k = 3$ with $\mathbf{m} = (0, m/2, n/2, 1)$ for some $n$ such that $m \leq n \leq 2$ and

$$A_3 = \begin{pmatrix} a & 0 \\ 0 & a \end{pmatrix}, \qquad Q^1 = \begin{pmatrix} q & q \\ q & q \end{pmatrix}, \qquad Q^2 = \begin{pmatrix} u & u \\ u & u \end{pmatrix}.$$

With this choice of parameters, (2.22) becomes $F_N(Q_\varepsilon) \leq U(n,a) + \mathcal{R}$ where

$$\begin{aligned} U(n,a) = &-\beta^2(\theta(1) - \theta(q)) - n\beta^2(\theta(u) - \theta(q)) - m\beta^2 \theta(q) \\ &+ \beta^2(\xi'(1) - \xi'(u)) + a - 1 - \log a \\ &- \frac{1}{n} \log \frac{a - n\beta^2(\xi'(u) - \xi'(q))}{a} \\ &- \frac{1}{m} \log \frac{a - n\beta^2(\xi'(u) - \xi'(q)) - m\beta^2 \xi'(q)}{a - n(\xi'(u) - \xi'(q))}. \end{aligned}$$

First of all, if we take $n = 1$ and

$$a = a_0 = b + \beta^2(\xi'(u) - \xi'(q)) = \frac{1}{1-q} + \beta^2(\xi'(u) - \xi'(q))$$

we get $U(1, a_0) = 2\mathcal{P}(\beta)$ by comparing with (4.8). Furthermore,

$$\begin{aligned} d(u) &= \left. \frac{\partial U(n, a_0)}{\partial n} \right|_{n=1} \\ &= -\beta^2(\theta(u) - \theta(q)) + \left(1 - \frac{1}{m}\right) \beta^2(\xi'(u) - \xi'(q))(1 - q) \\ &\quad + \frac{1}{m} \frac{\beta^2(\xi'(u) - \xi'(q))(1-q)}{1 - m\beta^2 \xi'(q)(1-q)} - \log(1 + \beta^2(\xi'(u) - \xi'(q))(1-q)) \end{aligned}$$

and we will show that $d(u) < 0$ for $u > q$. In order to prove this, we will notice that $d(q) = 0$ and show that $d'(u) < 0$ for $u > q$. Indeed,

$$d'(u) = \beta^2 \xi''(u) \left( -u + \frac{\beta^2(\xi'(u) - \xi'(q))(1-q)^2}{1 + \beta^2(\xi'(u) - \xi'(q))(1-q)} + \frac{\beta^2 \xi'(q)(1-q)^2}{1 - m\beta^2 \xi'(q)(1-q)} \right).$$

The last term can be simplified using the fact that the first equation in (4.1) is equivalent to

$$\frac{\beta^2 \xi'(q)(1-q)^2}{1 - m\beta^2 \xi'(q)(1-q)} = q$$

and then some simple algebra gives that for $u > q$, $d'(u) < 0$ if and only if

$$\beta^2 \xi'(u) - \beta^2 \xi'(q) < \frac{1}{1-u} - \frac{1}{1-q}.$$

Since two sides are equal for $u = q$ it is enough to show that for $u > q$

$$\beta^2 \xi''(u) < \frac{1}{(1-u)^2}.$$

If we consider the function

$$F(u) = \beta^2 \xi'(u) + \frac{1}{m}\left(\frac{1}{1-q+mq} - \frac{1}{1-q+m(q-u)}\right)$$

then it is easy to check that the critical point conditions (4.1) are equivalent to

$$F(0) = F(q) = 0, \qquad \int_0^q F(u)\,du = 0.$$

Therefore,

$$F'(u) = \beta^2 \xi''(u) - \frac{1}{(1-q+m(q-u))^2} = 0$$

has at least two solutions on $[0, q]$. But this equation can be rewritten as

$$(\beta^2 \xi''(u))^{-1/2} = 1 - q + m(q-u)$$

and since the right-hand side is convex in $u$, this equation has at most two solutions. Therefore, for $u \geq q$,

$$(\beta^2 \xi''(u))^{-1/2} \geq 1 - q + m(q-u) \geq 1 - u,$$

which implies

$$\beta^2 \xi''(u) < \frac{1}{\hat{x}^2(u)} = \frac{1}{(1-u)^2}.$$

This finally proves that for $u > q$ by fluctuating parameter $n$ around 1 we get that $F_N(Q_\varepsilon) < 2\mathcal{P}(\beta) + \mathcal{R}$ and this completes the proof of Theorem 4. □

**5. Proof of Theorem 2.** Let us first recall the main steps in the proof of Theorem 1 for $Q$ such that $\min r_j \geq 1$ when (1.12) becomes $F_N(Q_\varepsilon) \leq n\mathcal{P}(\beta) + \mathcal{R}$. The remark following the proof of Theorem 1 states that we use interpolation (2.3) with $k = 1$, Gaussian random vectors $(z_{0,i}^1, \ldots, z_{0,i}^n)$ with covariance $\Delta_0 = \beta^2 Q^1$ where $Q^1 = Q - \beta^{-1} I$ is defined in (3.11) and Gaussian random vectors $(z_{1,i}^1, \ldots, z_{1,i}^n)$ with covariance $\Delta_1 = \beta I$. Then for the function $\varphi(t)$ defined in (2.6), Theorem 5 gives

$$\varphi'(t) = -\frac{n}{2}\beta^2(\theta(1) - \theta(q)) - R(t) + \mathcal{R}, \tag{5.1}$$

where the remainder $R(t)$ was defined in (2.9). For $\xi(q) = q^2/2$

$$R(t) = \tfrac{1}{4}\beta^2 \sum_{j,j' \leq n} \mathbb{E}\langle (R^{j,j'} - q^1_{j,j'})^2 \rangle_t. \tag{5.2}$$

In order to show that the upper bound of Theorem 1 is sharp we need to control $R(t)$ and also prove that the bound on $\varphi(0)$ given in Lemma 1 is sharp. Since

$$\mathbb{E}\left(\sum_{j \leq n}\sum_{i \leq N} \sigma_i^j z_{1,i}^j\right)^2 = Nn\beta$$

we have

$$\varphi(0) = \frac{1}{N}\mathbb{E}\log \mathbb{E}_1 \int_{Q_\varepsilon} \exp \sum_{j \leq n}\sum_{i \leq N} \sigma_i^j(z_{0,i}^j + z_{1,i}^j)\, d\lambda_N^n$$

$$= \frac{n\beta}{2} + \frac{1}{N}\mathbb{E}\log \int_{Q_\varepsilon} \exp \sum_{j \leq n}\sum_{i \leq N} \sigma_i^j z_{0,i}^j\, d\lambda_N^n.$$

On the other hand, (2.11) gives, using $A_1$ as a parameter as in (2.22), that is, writing $A_2 = A_1 + \Delta_1$,

$$\varphi(0) \leq \frac{1}{2}(\operatorname{Tr}(Q\Delta_1) + \operatorname{Tr}(A_1 Q) + \operatorname{Tr}(A_1^{-1}\Delta_0) - n - \log|A_1|) + \mathcal{R}$$
$$= \frac{n\beta}{2} + \frac{1}{2}(\operatorname{Tr}(A_1 Q) + \operatorname{Tr}(A_1^{-1}\Delta_0) - n - \log|A_1|) + \mathcal{R}. \tag{5.3}$$

Therefore, in order to show that this bound is sharp, we need to prove the following.

LEMMA 4. *We have*

$$\lim_{\varepsilon \to 0}\lim_{N \to \infty} \frac{1}{N}\mathbb{E}\log \int_{Q_\varepsilon} \exp \sum_{j \leq n}(\boldsymbol{\sigma}^j, \mathbf{z}^j)\, d\lambda_N^n$$

$$= \inf_A \frac{1}{2}(\operatorname{Tr}(AQ) + \operatorname{Tr}(A^{-1}\Delta_0) - n - \log|A|),$$

*where* $\mathbf{z}^j = (z_{0,1}^j, \ldots, z_{0,N}^j)$.

The statement of the lemma will be proved for any $Q$ and $\Delta_0$, but for specific choices in the above interpolation it is easy to check that the infimum of (5.3) is achieved on

(5.4) $\quad A_1 = \beta I$ and $\displaystyle\lim_{\varepsilon \to 0} \lim_{N \to \infty} \varphi(0) = \frac{n}{2}(3\beta - 2 - \log \beta).$

PROOF OF LEMMA 4. The upper bound is given by (5.3). To prove the lower bound we will first replace integration over the sphere by a Gaussian integral constrained to the small neighborhood of the sphere. Let $\nu$ be a standard Gaussian distribution on $\mathbb{R}^N$. If we write $\boldsymbol{\rho} \in \mathbb{R}^N$ as $\boldsymbol{\rho} = s\boldsymbol{\sigma}$ for $\boldsymbol{\sigma} \in S_N$ then $\boldsymbol{\sigma}$ and $s$ are independent and the law of $\boldsymbol{\sigma}$ is $\lambda_N$. Denote by $\gamma$ the law of $s$. Clearly, $\gamma([1-\varepsilon, 1+\varepsilon]) \to 1$ for any $\varepsilon > 0$. Let us consider a set

(5.5) $\quad \Omega_\varepsilon = \{(\boldsymbol{\rho}^1, \ldots, \boldsymbol{\rho}^n) : \forall j \leq n, s_j \in [1-\varepsilon, 1+\varepsilon], (\boldsymbol{\sigma}^1, \ldots, \boldsymbol{\sigma}^n) \in Q_\varepsilon\}.$

First of all, for $(\boldsymbol{\rho}^1, \ldots, \boldsymbol{\rho}^n) \in \Omega_\varepsilon$,

$$\left| \sum_{j \leq n} ((\boldsymbol{\rho}^j, \mathbf{z}^j) - (\boldsymbol{\sigma}^j, \mathbf{z}^j)) \right| \leq \varepsilon \sqrt{N} \sum_{j \leq n} \|\mathbf{z}^j\|.$$

Therefore,

$$\int_{Q_\varepsilon} \exp \sum_{j \leq n} (\boldsymbol{\sigma}^j, \mathbf{z}^j) \, d\lambda_N^n$$

$$\geq \gamma([1-\varepsilon, 1+\varepsilon])^{-n} \exp\left(-\varepsilon \sqrt{N} \sum_{j \leq n} \|\mathbf{z}^j\|\right)$$

$$\times \int_{[1-\varepsilon, 1+\varepsilon]^n} \int_{Q_\varepsilon} \exp \sum_{j \leq n} (s_j \boldsymbol{\sigma}^j, \mathbf{z}^j)$$

$$d\lambda_N(\boldsymbol{\sigma}^1) \cdots d\lambda_N(\boldsymbol{\sigma}^n) \, d\gamma(s_1) \cdots d\gamma(s_n)$$

$$= \gamma([1-\varepsilon, 1+\varepsilon])^{-n} \exp\left(-\varepsilon \sqrt{N} \sum_{j \leq n} \|\mathbf{z}^j\|\right)$$

$$\times \int_{\Omega_\varepsilon} \exp \sum_{j \leq n} (\boldsymbol{\rho}^j, \mathbf{z}^j) \, d\nu^n(\boldsymbol{\rho}^1, \ldots, \boldsymbol{\rho}^n)$$

and

$$\frac{1}{N} \mathbb{E} \log \int_{Q_\varepsilon} \exp \sum_{j \leq n} (\boldsymbol{\sigma}^j, \mathbf{z}^j) \, d\lambda_N^n \geq \frac{1}{N} \mathbb{E} \log \int_{\Omega_\varepsilon} \exp \sum_{j \leq n} (\boldsymbol{\rho}^j, \mathbf{z}^j) \, d\nu^n - \mathcal{R}.$$

We can now replace the set $\Omega_\varepsilon$ by the set defined in terms of the overlaps of $\boldsymbol{\rho}^j$, $j \leq n$. Let us consider a set

$$\Omega_\delta = \{(\boldsymbol{\rho}^1, \ldots, \boldsymbol{\rho}^n) : R(\boldsymbol{\rho}^j, \boldsymbol{\rho}^{j'}) \in [q_{j,j'} - \delta, q_{j,j'} + \delta] \text{ for } j, j' \leq n\}.$$

Then, clearly, if we choose $\delta = \delta(\varepsilon)$ small enough then $\Omega_\delta \subseteq \Omega_\varepsilon$ and, therefore,

$$\frac{1}{N}\mathbb{E}\log \int_{Q_\varepsilon} \exp \sum_{j \leq n}(\boldsymbol{\sigma}^j, \mathbf{z}^j)\, d\lambda_N^n \geq \frac{1}{N}\mathbb{E}\log \int_{\Omega_\delta} \exp \sum_{j \leq n}(\boldsymbol{\rho}^j, \mathbf{z}^j)\, d\nu^n - \mathcal{R}.$$

Thus, we replaced integration over the sphere by a Gaussian integral. Given a set $V \subseteq (\mathbb{R}^N)^n$ and a symmetric $n \times n$ matrix $A$ such that $A+I$ is positive definite we define

$$\Phi_A(V) := \frac{1}{N}\mathbb{E}\log \int_V \exp\left(\sum_{j \leq n}(\boldsymbol{\rho}^j, \mathbf{z}^j) - \frac{1}{2}\sum_{j,j' \leq n} a_{j,j'}(\boldsymbol{\rho}^j, \boldsymbol{\rho}^{j'})\right) d\nu^n.$$

First of all,

(5.6) $\quad \Phi_A((\mathbb{R}^N)^n) = F(A) := \frac{1}{2}(\text{Tr}((A+I)^{-1}\Delta_0) - \log|A+I|),$

which is easy to show by decoupling the coordinates and proceeding as in (2.17). Using the definition of $\Omega_\delta$, we have

$$\frac{1}{N}\mathbb{E}\log \int_{\Omega_\delta} \exp \sum_{j \leq n}(\boldsymbol{\rho}^j, \mathbf{z}^j)\, d\nu^n \leq \frac{1}{2}\text{Tr}(AQ) + \Phi_A((\mathbb{R}^N)^n) + \mathcal{R}$$

$$= \frac{1}{2}\text{Tr}(AQ) + F(A) + \mathcal{R},$$

which, of course, coincides with the upper bound by making the change of variable $A \to A+I$. Let $A_0$ be a matrix that minimizes $\text{Tr}(AQ)/2 + F(A)$. Such $A_0$ exists and is unique since $F(A)$ is convex in $A$ by Hölder's inequality, the set $A+I > 0$ is convex and $F(A) \to +\infty$ when $|A+I| \to 0$. The critical point condition for $A_0$ is

(5.7) $\quad \left.\frac{\partial}{\partial x}\left(\frac{1}{2}\text{Tr}((A_0 + xB)Q) + F(A_0 + xB)\right)\right|_{x=0} = 0$

for any symmetric matrix $B$. Let us consider the sets

(5.8) $\quad V_{j,j'}^+ = \{R(\boldsymbol{\rho}^j, \boldsymbol{\rho}^{j'}) \geq q_{j,j'} + \delta\}, \qquad V_{j,j'}^- = \{R(\boldsymbol{\rho}^j, \boldsymbol{\rho}^{j'}) \leq q_{j,j'} - \delta\}$

for $j, j' \leq n$, so that

(5.9) $\quad (\mathbb{R}^N)^n = \Omega_\delta \cup \left(\bigcup_{j,j' \leq n} V_{j,j'}^+\right) \cup \left(\bigcup_{j,j' \leq n} V_{j,j'}^-\right).$

We will show below that if $V$ is any one of the sets (5.8) then for some $c > 0$

(5.10) $\quad \Phi_{A_0}(V) \leq \Phi_{A_0}((\mathbb{R}^N)^n) - c = F(A_0) - c.$

This implies that

(5.11) $\quad \lim_{N \to \infty} \Phi_{A_0}(\Omega_\delta) = \lim_{N \to \infty} \Phi_{A_0}((\mathbb{R}^N)^n) = F(A_0).$

Suppose this is not the case. Then for some positive $c > 0$ and $N$ large enough we would have

(5.12) $$\Phi_{A_0}(\Omega_\delta) \leq \Phi_{A_0}((\mathbb{R}^N)^n) - c.$$

Then by concentration of measure (5.10) and (5.12) would imply that with probability exponentially close to one for $V$ equal to $\Omega_\delta$ or one of the sets in (5.8)

$$\int_V \exp\left(\sum_{j\leq n}(\boldsymbol{\rho}^j, \mathbf{z}^j) - \tfrac{1}{2}\sum_{j,j'\leq n} a_{j,j'}(\boldsymbol{\rho}^j, \boldsymbol{\rho}^{j'})\right) d\nu^n$$
$$\leq \exp(-N/L) \int_{(\mathbb{R}^N)^n} \exp\left(\sum_{j\leq n}(\boldsymbol{\rho}^j, \mathbf{z}^j) - \tfrac{1}{2}\sum_{j,j'\leq n} a_{j,j'}(\boldsymbol{\rho}^j, \boldsymbol{\rho}^{j'})\right) d\nu^n.$$

Adding up all these inequalities and using (5.9) we arrive at a contradiction, proving (5.11), which in turn implies that

$$\frac{1}{N}\mathbb{E}\log \int_{\Omega_\delta} \exp \sum_{j\leq n}(\boldsymbol{\rho}^j, \mathbf{z}^j)\, d\nu^n = \frac{1}{2}\operatorname{Tr}(A_0 Q) + \Phi_{A_0}(\Omega_\delta) + \mathcal{R}$$
$$= \frac{1}{2}\operatorname{Tr}(A_0 Q) + F(A_0) + \mathcal{R}$$

and this will finishes the proof of the lemma. It remains only to prove (5.10). The proof is the same for all $V_{j,j'}^+$ and $V_{j,j'}^-$ so we will only consider the case $V = V_{j,j'}^-$ for $j \neq j'$. On the set $V_{j,j'}^-$ for any $x \geq 0$,

$$xR_{j,j'} \leq x(q_{j,j'} - \delta)$$

and, therefore, if we consider a matrix $B$ with two nonzero entries $b_{j,j'} = b_{j',j} = 1$, then

$$\Phi_{A_0}(V_{j,j'}^-) \leq x(q_{j,j'} - \delta) + \Phi_{A_0+xB}(V_{j,j'}^-) \leq x(q_{j,j'} - \delta) + F(A_0 + xB)$$
$$= U(x)$$
$$:= -x\delta - \tfrac{1}{2}\operatorname{Tr}(A_0 Q) + \tfrac{1}{2}\operatorname{Tr}((A_0 + xB)Q) + F(A_0 + xB).$$

We have $U(0) = F(A_0)$ and using (5.7), $U'(0) = -\delta$. Therefore, by slightly increasing $x$ we get (5.10) and this completes the proof of lemma. □

Next, we need to learn how to control the remainder term (5.2). The general approach is the same as in the proof of the Parisi formula in [7] (and in [8] or [4]), but the a priori estimates are now performed on $2n$ coupled copies of the system and will be similar in spirit to Theorem 1. Let us consider two $2n \times 2n$ block matrices

$$G^1 = \begin{pmatrix} Q^1 & Q^1 \\ Q^1 & Q^1 \end{pmatrix} \quad \text{and} \quad G^2 = \begin{pmatrix} Q^2 & Q^1 \\ Q^1 & Q^2 \end{pmatrix}$$

and let $(z_{0,i}^1, \ldots, z_{0,i}^{2n})$ and $(z_{1,i}^1, \ldots, z_{1,i}^{2n})$ be independent Gaussian vectors with covariances

$$\Delta_0 = \beta^2 G^1 \quad \text{and} \quad \Delta_1 = \beta^2(G^2 - G^1) = \beta^2(1-q)I = \beta I$$

correspondingly, independent for $i \leq N$. Consider

(5.13) $\quad H_t(\boldsymbol{\sigma}^1, \ldots, \boldsymbol{\sigma}^{2n}) = \sqrt{t} \sum_{j \leq 2n} \beta H_N(\boldsymbol{\sigma}^j) + \sqrt{1-t} \sum_{j \leq 2n} \sum_{i \leq N} \sigma_i^j (z_{0,i}^j + z_{1,i}^j)$

and

(5.14) $\quad h_t(\boldsymbol{\sigma}^1, \ldots, \boldsymbol{\sigma}^{2n}) = \sqrt{t} \sum_{j \leq 2n} \beta H_N(\boldsymbol{\sigma}^j) + \sqrt{1-t} \sum_{j \leq 2n} \sum_{i \leq N} \sigma_i^j z_{0,i}^j.$

The Gibbs average $\langle \cdot \rangle_t$ in (5.2) can be rewritten as

$$\langle f(\boldsymbol{\sigma}^1, \ldots, \boldsymbol{\sigma}^{2n}) \rangle_t = \int_{Q_\varepsilon \times Q_\varepsilon} f \exp h_t \, d\lambda_N^{2n} \Big/ \int_{Q_\varepsilon \times Q_\varepsilon} \exp h_t \, d\lambda_N^{2n}$$

and the remainder (5.2) can be rewritten as

(5.15) $\quad R(t) = \tfrac{1}{4}\beta^2 \sum_{1 \leq j \leq n} \sum_{n < j' \leq 2n} \mathbb{E}\langle (R_{j,j'} - g_{j,j'}^1)^2 \rangle_t.$

Given a set $A \subseteq S_N^{2n}$, consider

(5.16)
$$\Psi(A, t) = \frac{1}{N} \mathbb{E} \log \mathbb{E}_1 \int_A \exp H_t \, d\lambda_N^{2n}$$
$$= n\beta^2(1-q) + \frac{1}{N} \mathbb{E} \log \int_A \exp h_t \, d\lambda_N^{2n}.$$

Consider $n \times n$ matrix $E = (e_{j,j'})$ with elements $e_{j,j'} \in [-1, 1]$ such that

$$U = \begin{pmatrix} Q^2 & E \\ E & Q^2 \end{pmatrix}$$

is a nonnegative definite matrix and consider the set

(5.17) $\quad U_\varepsilon = \{(\boldsymbol{\sigma}^1, \ldots, \boldsymbol{\sigma}^{2n}) \in S_N^{2n} : R_{j,j'} \in [u_{j,j'} - \varepsilon, u_{j,j'} + \varepsilon] \text{ for } j, j' \leq 2n\}.$

Let us define

(5.18)
$$\psi(t) = \frac{n}{2}(3\beta - 2 - \log \beta - \beta^2 t(\theta(1) - \theta(q)))$$
$$= \frac{n}{2}(3\beta - 2 - \log \beta - t(\beta - 1/2)).$$

Theorem 2 will then follow from the following a priori estimate.

THEOREM 6. *For any $t_0 < 1$ and for any $t \leq t_0$, for $N$ large enough,*

$$\Psi(U_\varepsilon, t) \leq 2\psi(t) - \frac{1}{K} \sum_{j,j' \leq n} (e_{j,j'} - q^1_{j,j'})^2 + \mathcal{R} \tag{5.19}$$

*where a constant $K$ does not depend on $N, t$ and $U$.*

Let us first show how Theorem 6 implies Theorem 2.

PROOF OF THEOREM 2. First of all,

$$\psi(1) = \frac{n}{2}\left(2\beta - \frac{3}{2} - \log\beta\right) = n\mathcal{P}(\beta). \tag{5.20}$$

Let us take $K$ as in (5.19) and for $\varepsilon_1 > 0$ let us define a set

$$\mathcal{V} = \left\{ E = (e_{j,j'}) \in [-1,1]^{n^2} : \sum_{j,j' \leq n} (e_{j,j'} - q^1_{j,j'})^2 \geq 2K(\psi(t) - \varphi(t)) + 2K\varepsilon_1 \right\}.$$

For any $E \in \mathcal{V}$, (5.19) implies that

$$\Psi(U_\varepsilon, t) \leq 2\psi(t) - \frac{1}{K}(2K(\psi(t) - \varphi(t)) + 2K\varepsilon_1) + \mathcal{R} \leq 2\varphi(t) - \varepsilon_1 \tag{5.21}$$

for large enough $N$ and small enough $\varepsilon$. Everywhere below let $L$ denote a constant that might depend on $\varepsilon_1$ and $K$ denote a constant independent of $\varepsilon_1$. Since

$$2\varphi(t) = \Psi(Q_\varepsilon \times Q_\varepsilon, t),$$

(5.21) and concentration of measure imply that

$$\mathbb{E}\langle I(U_\varepsilon)\rangle_t \leq L\exp(-N/L) \tag{5.22}$$

and the constant $L$ here does not depend on $E$. Let us consider the set

$$\Omega = \left\{ \sum_{1 \leq j \leq n} \sum_{n < j' \leq 2n} (R_{j,j'} - g^1_{j,j'})^2 \geq 2K(\psi(t) - \varphi(t)) + 2K\varepsilon_1 \right\} \cap (Q_\varepsilon \times Q_\varepsilon).$$

We can choose a sequence $E_i \in \mathcal{V}$ for $i \leq K(\varepsilon)$ and large enough $K(\varepsilon)$ such that

$$\Omega \subseteq \bigcup_{i \leq K(\varepsilon)} U_\varepsilon(E_i),$$

where we made the dependence of $U_\varepsilon$ on $E$ explicit. Inequality (5.22) implies that

$$\mathbb{E}\langle I(\Omega)\rangle_t \leq LK(\varepsilon)\exp(-N/L). \tag{5.23}$$

Using the definition of $\psi$ in (5.18) and (5.1),

(5.24) $$(\psi(t) - \varphi(t))' = R(t) + \mathcal{R}.$$

On the complement $\Omega^c$ of $\Omega$ we have

$$\sum_{1 \leq j \leq n} \sum_{n < j' \leq 2n} (R_{j,j'} - g^1_{j,j'})^2 \leq 2K(\psi(t) - \varphi(t)) + 2K\varepsilon_1,$$

and (5.23) implies that

$$R(t) \leq K((\psi(t) - \varphi(t)) + \varepsilon_1 + \mathbb{E}\langle I(\Omega)\rangle_t)$$
$$\leq K(\psi(t) - \varphi(t)) + K\varepsilon_1 + LK(\varepsilon)\exp(-N/L).$$

The relation (5.24) now implies that

$$(\psi(t) - \varphi(t))' \leq K(\psi(t) - \varphi(t)) + K\varepsilon_1 + LK(\varepsilon)\exp(-N/L) + \mathcal{R}.$$

Since by Lemma 4 we have $\lim_{\varepsilon \to 0} \lim_{N \to \infty} \varphi(0) = \psi(0)$, solving this differential inequality and then letting $N \to \infty, \varepsilon \to 0$ and $\varepsilon_1 \to 0$ implies that

$$\lim_{\varepsilon \to 0} \lim_{N \to \infty} \varphi(t) = \psi(t) \quad \text{for } t \leq t_0.$$

Since the derivatives $\psi'(t)$ and $\varphi'(t)$ are both bounded we have

$$\limsup_{\varepsilon \to 0} \limsup_{N \to \infty} |\varphi(1) - \psi(1)| \leq K(1 - t_0).$$

Using (5.20),

$$\limsup_{\varepsilon \to 0} \limsup_{N \to \infty} |\varphi(1) - n\mathcal{P}(\beta)| \leq K(1 - t_0)$$

and letting $t_0 \to 1$ completes the proof of Theorem 2. $\square$

PROOF OF THEOREM 6. Let $U^2 = U$ and let $U^1$ be a symmetric nonnegative definite matrix that will be specified later such that $U^2 - U^1$ is also nonnegative definite. Let

$$\Delta'_0 = \beta^2 U^1 \quad \text{and} \quad \Delta'_1 = \beta^2(U^2 - U^1).$$

Let $(y^1_{0,i}, \ldots, y^{2n}_{0,i})$ and $(y^1_{1,i}, \ldots, y^{2n}_{1,i})$ be Gaussian vectors with covariances $\Delta'_0$ and $\Delta'_1$ correspondingly, independent for $i \leq N$. For $0 \leq s \leq 1$, consider

(5.25)
$$H_s(\boldsymbol{\sigma}^1, \ldots, \boldsymbol{\sigma}^{2n}) = \sqrt{st} \sum_{j \leq 2n} \beta H_N(\boldsymbol{\sigma}^j)$$
$$+ \sqrt{1-s}\sqrt{t} \sum_{j \leq 2n} \sum_{i \leq N} \sigma^j_i (y^j_{0,i} + y^j_{1,i})$$
$$+ \sqrt{1-t} \sum_{j \leq 2n} \sum_{i \leq N} \sigma^j_i (z^j_{0,i} + z^j_{1,i})$$

and let

$$\phi(s) = \frac{1}{N}\mathbb{E}\log\mathbb{E}_1\int_{U_\varepsilon}\exp H_s\,d\lambda_N^{2n}.$$

By a straightforward computation as in Theorem 5 one can show that

$$\phi'(s) = -\frac{t\beta^2}{2}\sum_{j,j'\leq 2n}(\theta(u^2_{j,j'}) - \theta(u^1_{j,j'})) - R(s) + \mathcal{R}$$

where the remainder $R(s) \geq 0$. The analogue of the first line in (2.7) is not present here because $U^2 = U$, that is, the covariance parameters match the constraints on the overlaps. Therefore,

$$(5.26)\quad \Psi(U_\varepsilon, t) = \phi(1) \leq \phi(0) - \frac{t\beta^2}{2}\sum_{j,j'\leq 2n}(\theta(u^2_{j,j'}) - \theta(u^1_{j,j'})) + \mathcal{R}.$$

Note that

$$H_0(\boldsymbol{\sigma}^1,\ldots,\boldsymbol{\sigma}^{2n}) = \sum_{j\leq 2n}\sum_{i\leq N}\sigma_i^j(x_{0,i}^j + x_{1,i}^j)$$

where $x_{l,i}^j = \sqrt{t}y_{l,i}^j + \sqrt{1-t}z_{l,i}^j$ and, therefore, the vectors $(x_{0,i}^1,\ldots,x_{0,i}^{2n})$ and $(x_{1,i}^1,\ldots,x_{1,i}^{2n})$ have covariances

$$\Delta_0^t = (1-t)\Delta_0 + t\Delta_0' \quad\text{and}\quad \Delta_1^t = (1-t)\Delta_1 + t\Delta_1'.$$

Lemma 1 now implies that for any $2n \times 2n$ symmetric positive definite matrix $A$ such that $A_1 = A - \Delta_1^t$ is also positive definite we have

$$2\phi(0) \leq \text{Tr}(\Delta_1^t U) + \text{Tr}(A_1 U) - 2n + \text{Tr}(A_1^{-1}\Delta_0^t) - \log|A_1| + \mathcal{R}$$

and combining with (5.26)

$$(5.27)\quad \begin{aligned}2\Psi(U_\varepsilon,t) \leq &\text{Tr}(\Delta_1^t U) + \text{Tr}(A_1 U) - 2n + \text{Tr}(A_1^{-1}\Delta_0^t) - \log|A_1|\\&- t\beta^2\sum_{j,j'\leq 2n}(\theta(u^2_{j,j'}) - \theta(u^1_{j,j'})) + \mathcal{R}.\end{aligned}$$

We proceed by a diagonalization procedure as in the proof of Theorem 1. Let $U^2 = O^T R^2 O$ for $R^2 = \text{Diag}(r_1^2,\ldots,r_{2n}^2)$ and orthogonal matrix $O$. Take $R^1 = \text{Diag}(r_1^1,\ldots,r_{2n}^1) > 0$ such that $R^2 - R^1$ is positive definite and let $U^1 = O^T R^1 O$. Take $B = \text{Diag}(b_1,\ldots,b_{2n})$ and let $A_1 = O^T BO$. Since $\Delta_1 = \beta^2(1-q)I$, we have $\Delta_1 = O^T \Delta_1 O$. Let $C$ be a matrix such that $\Delta_0 = \beta^2 O^T CO$. Then

$$\Delta_0^t = \beta^2 O^T((1-t)C + tR^1)O$$

and

$$\Delta_1^t = \beta^2 O^T((1-t)(1-q)I + t(R^2 - R^1))O.$$

Since for $\theta(q) = q^2/2$

$$\sum_{j,j' \leq 2n} \theta(u_{j,j'}^l) = \sum_{j \leq 2n} \theta(r_j^l),$$

the bound (5.27) becomes

$$2\Psi(U_\varepsilon, t) \leq \sum_{j \leq 2n} \Big( \beta^2(t(r_j^2 - r_j^1) + (1-t)(1-q))r_j^2 + b_j r_j^2 - 1$$

(5.28)
$$+ \frac{1}{b_j} \beta^2 (tr_j^1 + (1-t)c_{j,j}) - \log b_j - t\beta^2(\theta(r_j^2) - \theta(r_j^1)) \Big)$$

$$+ \mathcal{R}.$$

We will take

(5.29)
$$b_j = \frac{1}{r_j^2 - r_j^1}$$

and the bound (5.28) becomes

$$2\Psi(U_\varepsilon, t) \leq \sum_{j \leq 2n} \Big( \frac{1}{2}\beta^2 t((r_j^2)^2 - (r_j^1)^2) + \frac{r_j^2}{r_j^2 - r_j^1} - 1$$

(5.30)
$$+ \beta^2(1-t)((1-q)r_j^2 + (r_j^2 - r_j^1)c_{j,j}) + \log(r_j^2 - r_j^1) \Big)$$

$$+ \mathcal{R}.$$

Let

(5.31)
$$r_j^1 = \begin{cases} r_j^2 - \frac{1}{\beta}, & \text{if } r_j^2 \geq \frac{1}{\beta}, \\ 0, & \text{if } r_j^2 < \frac{1}{\beta}. \end{cases}$$

In the first case the $j$th term in the sum (5.30) becomes

$$f_1(r_j^2) = \beta t r_j^2 - \frac{t}{2} + \beta r_j^2 - 1 + \beta^2(1-t)(1-q)r_j^2 + \beta(1-t)c_{j,j} - \log \beta$$

$$= \beta t r_j^2 - \frac{t}{2} + \beta r_j^2 - 1 + \beta(1-t)r_j^2 + \beta(1-t)c_{j,j} - \log \beta$$

since $1 - q = 1/\beta$, and in the second case it becomes

$$f_2(r_j^2) = \frac{\beta^2 t}{2}(r_j^2)^2 + \beta(1-t)r_j^2 + \beta^2(1-t)r_j^2 c_{j,j} + \log r_j^2.$$

It is easy to check that

$$f_1(\beta^{-1}) = f_2(\beta^{-1}) \quad \text{and} \quad f_2'(\beta^{-1}) - f_1'(\beta^{-1}) = \beta^2(1-t)c_{j,j} \geq 0$$

since $C$ is nonnegative definite. Also,

$$f_2''(r_j^2) = \beta^2 t - \frac{1}{(r_j^2)^2} \leq -\beta^2(1-t) \leq 0 \qquad \text{for } r_j^2 \leq \frac{1}{\beta}$$

and $f_1''(r_j^2) = 0$ and, therefore,

$$(5.32) \qquad f_2(r_j^2) < f_1(r_j^2) \qquad \text{for } r_j^2 < \frac{1}{\beta}.$$

Therefore,

$$(5.33) \qquad 2\Psi(U_\varepsilon, t) \leq \sum_{j \leq 2n} f_1(r_j^2) + \mathcal{R}$$

and this bound is achieved by taking $r_j^1 = r_j^2 - 1/\beta$ in (5.30). It is important to note that in this bound we no longer have to assume that $r_j^2 \geq 1/\beta$ and even though for $r_j^2 < 1/\beta$ parameter $r_j^1 = r_j^2 - 1/\beta$ becomes negative, due to (5.32) we can simply treat (5.33) as a formula. Since

$$\sum_{j \leq 2n} r_j^2 = \operatorname{Tr}(R^2) = \operatorname{Tr}(U^2) = 2n \quad \text{and} \quad \sum_{j \leq 2n} c_{j,j} = \operatorname{Tr}(C) = \operatorname{Tr}(G^1) = 2nq$$

we have

$$\sum_{j \leq 2n} f_1(r_j^2) = 2n\beta t - nt + 2n\beta - 2n + 2n\beta(1-t) + 2n\beta(1-t)q - 2n\log\beta$$

$$= 2n(3\beta - 2 - \log\beta) - 2nt(\beta - 1/2) = 4\psi(t)$$

by comparing with (5.18). Therefore, (5.33) implies that

$$(5.34) \qquad \Psi(U_\varepsilon, t) \leq 2\psi(t) + \mathcal{R}.$$

However, we can improve upon (5.34) by slight fluctuations of $A_1$ in (5.27). Our choice of $r_j^1 = r_j^2 - 1/\beta$ and (5.29) imply that $b_j = \beta$, that is, $B = \beta I$ and, therefore, in (5.27) we have $A_1 = \beta I$. Also, $R^2 - R^1 = \beta^{-1}I$ and, hence, $U^1 = U^2 - \beta^{-1}I$. Let us take the derivative of (5.27) in $A_1$ at $A_1 = \beta I$. If $A_1$ was not constrained to be symmetric, we would have

$$\frac{\partial}{\partial A_1} \text{ r.h.s. of } (5.27)\Big|_{A_1 = \beta I}$$

$$= (U - (A_1^{-1})^T - (A_1^{-1}\Delta_0^t A_1^{-1})^T)\Big|_{A_1 = \beta I} = U^2 - \beta^{-1}I - \beta^{-2}\Delta_0^t$$

$$= U^1 - (1-t)G^1 - tU^1 = (1-t)(G^1 - U^1),$$

where we used some well-known formulas for derivatives of matrix determinants and inverses. However, since $A_1 = (a_{j,j'})$ is symmetric, $a_{j,j'} = a_{j',j}$, the derivative in the off-diagonal element will simply be doubled, that is, for $j \neq j'$,

$$d_{j,j'} = \frac{\partial}{\partial a_{j,j'}} \text{ r.h.s. of } (5.27)\Big|_{A_1=\beta I} = 2(1-t)(G^1 - U^1)_{j,j'}.$$

Since

$$G^1 - U^1 = (G^2 - \beta^{-1}I) - (U^2 - \beta^{-1}I) = G^2 - U^2$$
$$= \begin{pmatrix} 0 & Q^1 - E \\ Q^1 - E & 0 \end{pmatrix},$$

for $1 \leq j \leq n$ and $n < j' \leq 2n$ we have

$$d_{j,j'} = 2(1-t)(g^1_{j,j'} - u^1_{j,j'}) = 2(1-t)(q^1_{j,j'-n} - e_{j,j'-n}).$$

Also, since the second derivatives in $a_{j,j'}$ are bounded in some neighborhood of $A_1 = \beta I$, this implies that

$$(5.35) \qquad \Psi(U_\varepsilon, t) \leq 2\psi(t) - \frac{1}{K} \sum_{j,j' \leq n} (q^1_{j,j'} - e_{j,j'})^2 + \mathcal{R}$$

and this completes the proof of Theorem 6. $\square$

**Acknowledgments.** We would like to thank the referee for some very helpful comments and suggestions, in particular, for bringing to our attention the results in [2] (see remark at the end of Section 3). We would also like to thank Richard Dudley for some helpful comments.

Department of Mathematics  
Massachusetts Institute of Technology  
77 Massachusetts Avenue 2–181  
Cambridge, Massachusetts 02139  
USA  
E-mail: panchenk@math.mit.edu  

Equipe d'Analyse de l'Institut Mathématique  
Universite Paris 6  
4 Place Jussieu  
75230 Paris Cedex 05  
France  
and  
Department of Mathematics  
Ohio State University  
Columbus, Ohio 43210  
USA  
E-mail: spinglass@talagrand.net